\documentclass[11pt]{article}
\newcommand{\documentdate}{11 II 2026}

\usepackage{a4wide,latexsym,amsmath,varioref,graphicx,xcolor,amssymb,dsfont}

\title{A Fast Newton Method Under Local Lipschitz Smoothness} 

\author{	
	S. Gratton%
	\thanks{Universit\'e de Toulouse, INP, IRIT, Toulouse, France. Email:
		serge.gratton@enseeiht.fr. Work partially supported by 3IA Artificial and
		Natural Intelligence Toulouse Institute (ANITI), French "Investing for the Future
		- PIA3" program under the Grant agreement ANR-19-PI3A-0004"}, 
	~S. Jerad%
	\thanks{ Mathematical Institute, University of Oxford, United Kingdom. Email:
		sadok.jerad@maths.ox.ac.uk. Work started in Toulouse. This work was supported by the Hong Kong Innovation and Technology Commission (InnoHK Project CIMDA). }
	~and Ph. L. Toint%
	\thanks{NAXYS, University of Namur, Namur, Belgium. Email:
		philippe.toint@unamur.be. Partly supported by ANITI.}
}

\newcommand{\beqn}[1]{\begin{equation}\label{#1}}
	\newcommand{\eeqn}{\end{equation}}
\newcommand{\req}[1]{(\ref{#1})}
\newcommand{\ms}{\;\;\;\;}
\newcommand{\tim}[1]{\;\; \mbox{#1} \;\;}

\newtheorem{theorem}{Theorem}[section]
\newtheorem{lemma}[theorem]{Lemma}
\newtheorem{corollary}{Corollary}[section]
\newtheorem{assumption}{Assumption}
\newcommand{\numsection}[1]{\section{#1}\setcounter{equation}{0}}

\newcounter{algo}[section]
\renewcommand{\thealgo}{\thesection.\arabic{algo}}
\newcommand{\llem}[2]{\vspace{\baselineskip} 
	\noindent\framebox[\textwidth]{\parbox{0.95\textwidth}{
			\begin{lemma} \label{#1} \rm #2 \end{lemma} } } \vspace{\baselineskip} }
\newcommand{\assump}[2]{\vspace{\baselineskip} 
	\noindent\framebox[\textwidth]{\parbox{0.95\textwidth}{
			\begin{assumption} \label{#1} \rm #2 \end{assumption} } } \vspace{\baselineskip} }
\newcommand{\algo}[3]{{\refstepcounter{algo}
		{\small
			\begin{center}\begin{figure}[htbp]
					\framebox[\textwidth]{
						\parbox{0.95\textwidth} {\vspace{\topsep}
							{\bf Algorithm \thealgo : #2}\label{#1}\\
							\vspace*{-\topsep} \mbox{ }\\
							{#3} \vspace{\topsep} }}
		\end{figure}\end{center}}
	}
}
\newcommand{\bpr}{{\bf Proof.} \hspace{1.5mm}}
\newcommand{\epr}{\hfill $\Box$ \vspace*{1em}}
\newcommand{\proof}[1]{
	\begin{list}{}{
			\setlength{\topsep}{0.0pt}
			\setlength{\partopsep}{0.0pt}
			\setlength{\leftmargin}{0.025\textwidth}
			\setlength{\rightmargin}{0.5\leftmargin}
			\setlength{\labelwidth}{0.5\leftmargin}
			\setlength{\labelsep}{0.25\leftmargin}}
		\item \bpr #1 \epr \noindent
\end{list}}
\newcommand{\lthm}[2]{\vspace{\baselineskip} 
	\noindent\framebox[\textwidth]{\parbox{0.95\textwidth}{
			\begin{theorem} \label{#1} \rm #2 \end{theorem} } } \vspace{\baselineskip} }

\newcommand{\ii}[1]{\{ 1, \ldots, #1 \}}
\newcommand{\iiz}[1]{\{ 0, \ldots, #1 \}}
\newcommand{\iibe}[2]{\{ #1, \ldots, #2 \}}

\newcommand{\calO}{{\cal O}} 
\newcommand{\calS}{{\cal S}}

\newcommand{\calI}{{\cal I}}

\renewcommand{\Re}{\hbox{I\hskip -2pt R}}
\newcommand{\bigfrac}[2]{\frac{\displaystyle #1}{\displaystyle #2}}
\newcommand{\sfrac}[2]{{\scriptstyle \frac{#1}{#2}}}
\newcommand{\half}{\sfrac{1}{2}}
\newcommand{\eqdef}{\stackrel{\rm def}{=}}

\newcommand{\kap}[1]{\kappa_{\mbox{\tiny #1}}}

\newcommand{\al}[1]{{\footnotesize{\sf #1}}}
\newcommand{\tal}[1]{{\normalsize {\sf #1}}}

\newcommand{\sqrtsigkgk}{\sqrt{\sigma_k} \| g_k\|}

\counterwithin*{equation}{section}

\newcommand{\comment}[1]{}

\date{\documentdate}

\begin{document}

\maketitle

\begin{abstract}
	A new, fast second-order method is proposed that achieves the optimal
	$\mathcal{O}\left(|\log(\epsilon)|\epsilon^{-3/2}\right)$
	complexity to obtain first-order $\epsilon$-stationary points. 
	Crucially, this is deduced without assuming the
	standard global Lipschitz Hessian continuity condition, but only 
	using an appropriate local smoothness requirement. The algorithm 
	exploits Hessian information to compute a Newton step and a negative 
	curvature step when needed, in an approach similar to that of
	\cite{GratJerToin23bIMA}. 
	Inexact versions of the Newton step and negative curvature are
	proposed in order to reduce the cost of evaluating second-order information.
	Details are given of such an iterative implementation using Krylov subspaces.
	An extended algorithm for finding second-order
	critical points is also developed and its complexity is again shown
	to be within a log factor of the optimal one.
	Initial numerical experiments are discussed for both factorised
	and Krylov variants, which demonstrate the competitiveness of the
	proposed algorithm.  
\end{abstract}

Standard nonlinear optimization algorithms such as gradient
descent or Newton's method are widely used in a vast array of
applications. However their theoretical analysis
typically uses the somewhat restrictive assumptions that
gradients or higher derivatives are globally Lipschitz continuous
\cite{Cartis2022-wb}. Unfortunately, this assumption does not hold in
simple cases such as polynomials or exponentials. While it can be
argued that the assumption is only necessary on the ``tree of iterates'' 
(that is the union of all segments defined by algorithm's steps 
(see \cite[Notes p.47]{Cartis2022-wb} for instance), it remains 
desirable to relax it more significantly. Some progress 
in this direction has been made in \cite{MalitskyMischenko20,ZhouMAyang24}, 
where only local Lipschitz smoothness is assumed, but it only applies to 
the convex case. 

In the nonconvex case, motivated by the geometry of the
loss-landscape in modern deep learning problems and the success of
normalized gradient descent in this context, \cite{Zhang2020Why} proposed
a new smoothness condition for the analysis of this latter method, called
$(L_0,L_1)$  smoothness, which assumes that there exist constants $L_0
\geq 0$ and $L_1 > 0$ , such that,  
\beqn{L0L1smooth}
\|\nabla_x^2 f(x)\| \leq L_0 + L_1 \|\nabla_x^1 f(x)\|.
\eeqn 
This new assumption was shown to be more relevant in practice as it
covers univariate polynomial functions of all degrees, for instance.
We refer the reader to \cite{Zhang2020Why,ZhangJunFangWang20} for
additional motivation for the study of first-order method under 
\eqref{L0L1smooth} and for a complexity
analysis in both the deterministic and stochastic cases.
This proposal generated a new interest in the analysis of nonconvex
gradient descent. New classes of smoothness
and refinements of the initial analysis have been proposed in
\cite{KoloHendrStich23,ZhangJunFangWang20,LiQianTianRakJad23,ChenZhouLiangLu23}
both in the exact and stochastic case. Extensions to other classes of
problems or different update rules have been proposed, such as
adaptive-gradient methods \cite{Fawetal23,LiRakhlinJadb23}, variance
reduction \cite{Reisizadehetal23} and variational problems \cite{SunKaraRicht23}
to name just a few. Note also that \eqref{L0L1smooth} can be rephrased
as
\beqn{L0L1smoothrelax}
\|\nabla_x^1 f(x) - \nabla_x^1 f(y)\|
\leq (L_0 + L_1 \|\nabla_x^1 f(x)\|) , \tim{ if } \|x-y\| \leq \frac{1}{L_1},
\eeqn
to avoid the use of second-order derivatives (see
\cite{Reisizadehetal23,ZhangJunFangWang20} for further discussions on
the equivalence between \eqref{L0L1smooth} and \eqref{L0L1smoothrelax}). 

Extending the approach to second-order methods for optimization is
thus of interest, in order to enlarge the applicability of such
results beyond what is used for standard approaches
\cite{Cartis2022-wb}.  To the best of our knowledge, only
\cite{XieLiZhangDengGeYe24} has so far proposed an optimal
trust-region method under a Lipschitz smoothness assumptions similar
to \eqref{L0L1smoothrelax}. However, the proposed algorithm
unrealistically requires the knowledge of the problem's Lipschitz
constant, uses a trust region with a fixed radius, and computes an
exact solution of the trust-region subproblem, thus requiring intensive
numerical linear algebra to compute the step at each iteration.

{On the other hand, a vast array of work \cite{Royer2019,YaoXuRoosMahoWri22,Zhou25a,GratJerToin23bIMA} has been proposed in order to develop fast regularized Newton methods with a complexity guarantee close to $\calO\left(\epsilon^{-3/2}\right)$ evaluations to find an
	$\epsilon$-first order stationary point. Among this work, some required the knowledge of the criticality threshold \cite{Royer2019,YaoXuRoosMahoWri22} to achieve the bound  $\calO\left(\epsilon^{-3/2}\right)$. This requirement was bypassed in \cite{GratJerToin23bIMA} by alternating carefully  between regularized Newton step and  negative curvature and utilizing the current gradient norm to regularize the Newton system. However, this approach incurs an additional $|\log(\epsilon)|$ in the complexity. Recently, \cite{Zhou25a} combined the regularization proposed in \cite{GratJerToin23bIMA} with the capped-CG analysis developed in \cite{Royer2019} and reduced the extra cost of not knowing the criticality threshold $\epsilon$ to $|\log (\log(\frac{1}{\epsilon}))|$ which is better than the $|\log(\epsilon)|$ of \cite{GratJerToin23bIMA}.     }

In the present paper, inspired by the new Hessian Lipschitz
smoothness condition proposed in \cite{XieLiZhangDengGeYe24}, we
develop an adaptive Newton's method that requires at most
$\calO\left(|\log(\epsilon)|\epsilon^{-3/2}\right)$ evaluations to find an
$\epsilon$-first order stationary point. The met hod proposed is in
spirit of \cite{GratJerToin23bIMA} as it
alternates between Newton directions and directions of negative
curvature and also uses negative-curvature information to regularize
the Newton step. Beyond requiring significantly weaker smoothness,
it however differs from the proposal in this 
reference in crucial aspects, such as the power of regularization 
and the mechanism to accept or reject trial steps. The new method 
is fully adaptive and knowledge of the problem's geometry is 
not assumed. Two implementations 
are proposed, the first requiring exact negative curvature computation
and exact linear-system solves, and the second allowing more
inexactness and exploiting the structure of nested Krylov subspaces.   

The paper is organized as follows. Section~\ref{thealgo-s} starts by
stating our new Lipschitz smoothness condition, describes the general
algorithmic framework, compares it with closely related work on
second-order methods and states properties on the computed
step. Section~\ref{complexity-s} derives a bound on its worst-case
complexity for finding first-order critical
points. Section~\ref{so-teaser} details an algorithmic enhancement for
the search for a second-order critical point, with its associated
complexity proof discussed in Appendix~{\ref{app:so-theory}.} Section~\ref{stepcomproutine} describes
proposed procedures to compute the step and the incorporation
of preconditionning. The numerical behaviour of the proposed algorithms is 
considered in Section~\ref{sec:numerics} and some conclusions and perspectives 
are discussed in Section~\ref{sec:conclusions}.

\paragraph*{Notations}
Let $n \geq 1$. The symbol $\|\cdot\|$ denotes the Euclidean norm
for vectors in $\Re^n$ and its associated subordinate norm for matrices.
$\lambda_{\min}(M)$  and $\lambda_{\max}(M)$ denote the minimum and
maximum eigenvalues of a symmetric matrix $M$, while $I_n$ is the
identity matrix in $\Re^{n \times n}$.  For $x \in \Re$, we define
$[x]_{+} = \max(x,0)$. For two vectors $x,y \in \Re^n $,
$x^\intercal y$ denotes their inner product. The $i$-th column of $I_n$ is
denoted by $e_i$.

\numsection{Adaptive Newton with Negative Curvature using \\
	Local Smoothness}\label{thealgo-s}

We consider the problem of finding approximate first-order critical
points of the smooth unconstrained nonconvex optimization problem
\beqn{minf}
\min_{x \in \Re^n} f(x)
\eeqn
under the following set of assumptions.

\noindent
\textbf{AS.1} The function $f$ is two times continuously differentiable in $\Re^n$. 

\noindent
{\bf AS.2} There exists a constant $f_{\rm low}$ such that
$f(x) \geq f_{\rm low}$ for all $x\in \Re^n$.

\noindent
\textbf{AS.3} There exist constants $L_0\geq0$ and  $L_1 {\geq} 0$ and
$\delta {>} 0$, such that, if $\|x-y\| \leq \delta$, then
\beqn{LipHessian}
\|\nabla_x^2 f(y) - \nabla_x^2 f(x)\| \leq (L_0  + L_1 \|\nabla_x^1 f(x)\| ) \|x-y\|.  
\eeqn

\noindent
\textbf{AS.4}
There exists a constant $\kappa_B \geq 0$ such that
\[
\max(0,-\lambda_{\min}(\nabla_x^2 f(x))) \leq \kappa_B
\tim{for all} x\in \{ y \in \Re^n \mid f(y)\leq f(x_0)\}.
\]
AS.1 and AS.2 are standard assumptions when analyzing algorithms that
use second-order information
\cite{CarGouTo11b,BirgGardMartSantToin17,Cartis2022-wb}. AS.4 has also been considered for the analysis of other fast second-order methods, see \cite{GratJerToin23bIMA} for more discussion. AS.4 can be replaced by \eqref{L0L1smooth} if requested.

We now discuss our newly proposed AS.3. First note that  we
recover the standard Lipschitz Hessian 
smoothness requirement when $L_1 = 0$ and $\delta = \infty$. {Therefore,  based on the standard complexity lower bound established in \cite{Carmon2019a} for Lipschitz Hessian-smooth function, the best achievable complexity iteration  to reach an $\epsilon$-first order stationary point scales as $\mathcal{O}\left(\epsilon^{-3/2}\right)$. }  

As
mentioned in the introduction, Assumption  \eqref{LipHessian} has
recently been proposed for the study of second-order methods under
local smoothness, see \cite{XieLiZhangDengGeYe24} for instance. We now
state sufficient conditions on the function $f$ so that
\eqref{LipHessian} holds. 

\llem{suffconds}{Let $f$ be three times differentiable and suppose
	that there exists $M_0\geq 0$ and $M_1 > 0$ such that
	\beqn{thirdordersmoo}
	\|\nabla_x^3 f(x)[u]\|
	\leq (M_0 + M_1\|\nabla_x^1 f(x)\|) \|u\| \tim{ for all } u \in \mathbb{R}^n,
	\eeqn
	and that $f$ verifies \eqref{L0L1smoothrelax}. Then, there exits $(L_0,L_1,\delta)$ such that AS.3 holds.
}     
\proof{See \cite[Lemma~C.1]{XieLiZhangDengGeYe24}}

As a consequence,  it can be proved, using the equivalence
between \eqref{L0L1smoothrelax} and \eqref{L0L1smooth} for twice
differentiable functions, the last Lemma, that the class of functions satisfying  AS.3 
contains homogeneous or univariate polynomials of arbitrary degree and univariate exponentials, among others.  Indeed, 
their definition implies that derivatives of higher degree grow (at infinity) 
slower than derivatives of lower degree. 

We now detail the impact of this new condition on standard function
and gradient Hessian error bound.  

\llem{Lipbounds}{
	Suppose that AS.1 and AS.3 hold. Let $x$ and $s$ be in $\Re^n$ and
	$\|s\| \leq \delta$. Then,
	\beqn{Lipf}
	f(x+s) - f(x) - s^\intercal \nabla_x^1 f(x) - \frac{1}{2} s^\intercal \nabla_x^2 f(x) s^\intercal \leq \frac{\left(L_0  + L_1 \|\nabla_x^1 f(x)\| \right)}{6} \|s\|^3,
	\eeqn 
	and
	\beqn{Lipg}
	\| \nabla_x^1 f(x+s) - \nabla_x^1 f(x) - \nabla_x^2 f(x)s \| \leq \frac{\left(L_0  + L_1  \|\nabla_x^1 f(x)\| \right)}{2} \|s\|^2.
	\eeqn
}
\proof{
	The proof is an adaptation of the standard one for Lipschitz bounds
	\cite[Lemma~2.1]{CartGoulToin20b} taking the new condition \eqref{LipHessian} 
	into account. It is given in Appendix~{\ref{app:so-theory}} for the sake of completeness.
} 
The motivation for our algorithm is similar to that used in \cite{GratJerToin23bIMA} and we first describe its mechanism in very broad lines.
At each iteration, the idea is to minimize a "doubly regularized" quadratic model of the form
\beqn{reg-model-def}
m_k(s) \eqdef g_k^\intercal s + \half s^\intercal(H_k {+}  (\sqrtsigkgk + \mu_k)I_n)s
\eeqn
where {$g_k \eqdef \nabla_x^1 f(x_k)$}, $H_k \eqdef \nabla_x^2f(x_k)$, $\sigma_k$ is an adaptive, iteration dependent regularization parameter and $\mu_k$ and additional (hopefully small) regularization term depending on the type of step taken. The algorithm attempts to select a "small" $\mu_k$  and an associated Newton-like trial step using the regularized Hessian $H_k + (\sqrtsigkgk + \mu_k)I_n$. If this $\mu_k$ is too large or its associated trial step cannot be accepted (according to criteria detailed below), then a trial step along an approximate negative curvature direction is attempted.
The algorithm then proceeds, as is typical for adaptive regularization methods, by deciding if the trial step can be accepted, or has to be rejected.  The parameter $\sigma_k$ is then updated before starting a new iteration.

The computation an a-priori regularization and an associated trial step is clearly important, and 
can be organized in more than one way. Because we have two different implementations in mind (one uses matrix factorizations and the other an iterative Krylov approach, see Section~\ref{stepcomproutine}), we take for now a more abstract point of view of this computation  and  only assume, at this stage, that there exist a \al{Stepcomp}  procedure, which, given the current gradient, Hessian and a few algorithmic constants, produces both a regularization parameter $\mu_k$ and a trial step $s_k^{trial}$. 

This remains admittedly vague at this stage, but hopefully provides some intuition for the formal definition of the \tal{AN2CLS} algorithm \vpageref{AN2Cgen}. 

{\small 
	\algo{AN2Cgen}{  Adaptive Newton with Negative Curvature using Local Smoothness  (\tal{AN2CLS})}{
		\vspace*{-3mm}
		\begin{description}
			{\small
				\item[Step 0: Initialization: ] An initial point $x_0\in \Re^n$, a
				regularization 
				parameter $\sigma_0>0$ and a gradient accuracy threshold
				$\epsilon \in (0,1]$  are given, as well as the parameters
				\[
				\sigma_{\min}  > 0, \, \, \kappa_C \geq 1, \, \,   \kappa_\theta \geq 0, \, 0< \theta \leq 1, \, \vartheta  \geq 1,  
				\] 
				\[
				0< \gamma_1 < 1 < \gamma_2 \leq \gamma_3 
				\tim{and }
				0 < \eta_1 \leq \eta_2 < 1.  
				\]
				\vspace*{-2mm}
				Set $k=0$, define  
				\[\kap{slow} \eqdef (1+\kappa_\theta+\kappa_C)
				+ \sqrt{(1+\kappa_\theta+\kappa_C)^2 + \vartheta},
				\ms
				\kap{upnewt} \eqdef 3(1-\eta_2)   + 1 + \kappa_C   +\kappa_\theta\
				\]
				and set \texttt{REJECT = FALSE}.
				\item[Step 1: Check termination: ] If not already available, evaluate
				${g_k} \eqdef \nabla_x^1 f(x_k)$ and terminate if
				$\|g_k \| \leq \epsilon$.  Otherwise, evaluate $H_k \eqdef \nabla_x^2 f(x_k)$.
				
				\item [Step 2: Compute a trial step: ]
				Call the step computation procedure
				\beqn{Computproc}
				(s_k^{trial}, \mu_k) = \tal{Stepcomp}(g_k,H_k, \sigma_{k}, \kappa_C,\kappa_\theta, \theta)
				\eeqn
				
				\item[Step 3: Newton Step: ] If $\mu_k \leq \kappa_C\sqrtsigkgk$, 
				then evaluate $\nabla_x^1f(x_k+s_k^{trial})$.
				If 
				\beqn{condreject}
				\|\nabla_x^1f(x_k+s_k^{trial})\| >\frac{\|g_k\|}{2}   
				\tim{ and } 
				\|s_k^{trial}\| < \frac{1}{\sqrt{\sigma_{k}}\kap{slow}},
				\eeqn
				then set \texttt{REJECT = TRUE} and go to Step~6.
				Otherwise, set 
				\beqn{newtstep}
				s_k = s_k^{trial} \tim{ and } \kappa_k =\kap{upnewt}
				\eeqn
				and go to Step~5.
				\item[Step 4: Negative curvature step: ]  If $\mu_k > \kappa_C \sqrtsigkgk$, set
				\beqn{negcurvstep}
				s_k = s_k^{trial}  
				\tim{and}
				\kap{k} 
				\eqdef \frac{3}{2}\kappa_C^2 \theta^2 (1-\eta_2) + 1 + \frac{\kappa_C \mu_k}{\sqrt{\sigma_{k}}}.
				\eeqn
				\item[Step 5: Acceptance ratio computation and gradient test: ]
				Evaluate $f(x_k +s_k)$ and compute the acceptance ratio
				\beqn{rhokdef}
				\rho_k = \frac{f(x_k) - f(x_k +s_k)}{-(g_k^\intercal s_k + \frac{1}{2} s_k^\intercal H_k s_k)}.
				\eeqn
				If $\rho_k < \eta_1$ or
				\beqn{nextgradcond}
				\|\nabla_x^1 f(x_k + s_k)\| > \kap{k} \frac{\|g_k\|}{\epsilon},
				\eeqn
				set  \texttt{REJECT = TRUE}.
				\item[Step 6: Variables' update:]
				If \texttt{REJECT = FALSE}, set $x_{k+1}=x_k+s_k$, otherwise set $x_{k+1}=x_k$.
				\item[Step 7: Regularization parameter update:] 
				
				Set
				\begin{equation}\label{sigmakupdate}
					\sigma_{k+1} \in \left\{
					\begin{aligned}
						& [ \max\left(\sigma_{\min}, \gamma_1\sigma_k \right),\sigma_k] 
						&&\text{ if } \rho_k \geq \eta_2 
						&&\text{and } \texttt{REJECT=FALSE}, \\
						& [\sigma_k ,\gamma_2\sigma_k ]            
						&&\text{ if } \rho_k \in [\eta_1, \eta_2) &&\text{and } \texttt{REJECT=FALSE},\\
						& [\gamma_2 \sigma_{k}, \gamma_3 \sigma_{k}] 
						&&\text{ if } \rho_k < \eta_1 
						&&\text{and } \texttt{REJECT=TRUE}, 
					\end{aligned}
					\right.
				\end{equation}
				Increment $k$ by one, set \texttt{REJECT = FALSE} and go to Step~1.
			}
		\end{description}
	}
}

We have already mentioned that we wish to postpone the details of the procedure to compute the trial step to Section~\ref{stepcomproutine}, but, if we wish to analyze the method, we clearly need to be more specific, if not on the procedure, at least on the result of this procedure.  This is the object of the next assumption.  We therefore assume the following.

\assump{assump0}{
	The \tal{Stepcomp} procedure computes  a tentative regularization parameter $\mu_k$ and a trial step $s_k^{trial}$ satisfying the {following} conditions.\\
	If $\mu_k \leq \kappa_C \sqrtsigkgk$, then $s_k^{trial}$ is such that
	\beqn{posdefnewton}
	(s_k^{trial})^\intercal (H_k + \mu_kI_n) s_k^{trial} \geq 0,
	\eeqn
	\beqn{newtcond}
	\|r_k^{trial}\|
	\leq  \kappa_{\theta} \min\left(  \sqrtsigkgk \|s_k^{trial}\|,\|g_k\|\right)
	\eeqn
	and
	\beqn{orthcond}
	(r_k^{trial})^\intercal s_k^{trial} = 0,
	\eeqn
	where
	\beqn{rneig-def}
	r_k^{trial} \eqdef (H_k + (\sqrtsigkgk + \mu_k)I_n)s_k^{trial} + g_k.
	\eeqn
	
	Else if $\mu_k >  \kappa_C \sqrtsigkgk$, $s_k^{trial}$ is given by
	\beqn{sncurv-def}
	s_k^{trial} = \frac{\theta \kappa_C}{\sqrt{\sigma_k}}u_k
	\eeqn
	where the vector $u_k$  verifies
	\beqn{negcurvvector}
	g_k^\intercal u_k  \leq 0, \, \,  \|u_k \| = 1, \, \,
	u_k^\intercal H_k u_k \leq -\theta \mu_k
	\tim{and}
	u_k^\intercal H_k^2 u_k \leq \frac{\mu_k^2}{\theta^2}.
	\eeqn
}

\noindent
Further comments are in order at this point.
\begin{enumerate}
	\item  The reader familiar with the \al{AN2C} algorithm proposed in \cite{GratJerToin23bIMA} will notice the similarities, but also the differences between this algorithm and \al{AN2CLS}, the most important being the power of the regularization. In \al{AN2C} and inspired by the analysis of the convex case \cite{Mishchenko2023,DoikNest23},  the regularization has the form $\sqrt{\sigma_{k}\|g_k\|}$ whereas
	the regularization in \al{AN2CLS} is of the order of $\|g_k\|$ as
	is clear from \eqref{reg-model-def}. Note that choosing a
	regularization proportional to $\|g_k\|$ has already been considered in
	\cite{DoikMischNest22,doikov2023minimizing} for the convex case, with good initial
	numerical results and optimal complexity rates. In
	\cite{DoikMischNest22} a regularization proportional to $\|g_k\|$ was
	shown to be universal for a large class of convex functions ranging from those with Lipschitz continuous gradient to those with Lipschitz continuous
	third-order derivative. In \cite{doikov2023minimizing}, it was proved
	to be optimal  for a class of quasi-self concordant functions (see the
	above references for more details). \al{AN2CLS}
	generalizes this proposal to the nonconvex case with local smoothness only,
	while keeping the global rate of convergence close to optimal, as we show below. 
	Older proposals considering a Newton step regularized with a $\|g_k\|$ term in a linesearch setting may be found in \cite{UedaYama14,Polyak2007} but the derived complexity is suboptimal in both cases. 
	
	\item The new method also differs from that of \cite{GratJerToin23bIMA} in the mechanism used to accept or reject the trial step. Because the Lipschitz condition is now only local, one needs to explicitly check that the gradient at the trial point does not grow out of reasonable bounds (see \req{nextgradcond}) or too fast compared to the length of the trial step to be accepted (see \req{condreject}).  Large gradients at the trial point can only be accepted if the step is sufficiently large. 
	
	\item When \al{Stepcomp} returns a value of $\mu_k \leq \kappa_C\sqrt{\sigma_k}\|g_k\|$, the trial step is interpreted as "Newton step" because \req{posdefnewton}-\req{rneig-def} (which must hold in that case, see Assumption~0) imply that it (approximately) minimizes the "doubly regularized" model \req{reg-model-def}. Considering both conditions \eqref{newtcond} and \eqref{orthcond} has already been proposed in \cite{DussaMigOrb23} in order to develop scalable variants of cubic regularization. 
	
	\item In Step~1 of \al{AN2CLS}, the words "If not already available" are justified by the observation that, if iteration $k-1$ selected the Newton step in \req{newtstep} and accepted it, the value of 
	$g_k = \nabla_x^1f(x_k)= \nabla_x^1f(x_{k-1}+s_{k-1}^{trial})$ has been computed at Step~3 of iteration $k-1$.
	
	\item Condition \req{newtcond} in the residual of the "Newton step" is looser than that used in \al{AN2C}, since $\kappa_\theta$ is restricted to the interval $[0,1)$ in that algorithm. Also note that the \al{AN2CLS} algorithm has fewer hyperparameters than \al{AN2C}.
\end{enumerate}

Before proceeding with the analysis, we will introduce some useful notation.
Following well-established practice, we define
\[
\calS \eqdef \{ k \geq 0 \mid x_{k+1} = x_k+s_k \} 
\]
the set of indexes of ``successful iterations'', and
\[
\calS_k \eqdef \calS \cap \iiz{k},
\]
the set of indexes of successful iterations up to iteration $k$. We
further partition the iterations into two subsets depending on the nature of the step taken, according to
\[
\calI^{newt} \eqdef \{i\geq0 \, | \,\mu_i \leq \kappa_C \sqrt{\sigma_i} \|g_i\| \} 
\tim{ and }
\calI^{ncurv} \eqdef \{i\geq 0 \, | \, \mu_i > \kappa_C \sqrt{\sigma_i} \|g_i\| \},
\]
the first set containing the indices of the iterations where $s_k^{trial}$ is a Newton step and the second where it is a negative curvature step.
Moreover, considering \eqref{condreject}, we further subdivide the subset of $\calI^{newt}$ into three subsets as follows,
\begin{align*}
	\calI^{g\searrow} &
	\eqdef \left\{i\in \calI^{newt} \, \Big| \; \|\nabla_x^1f(x_i+s_i^{trial})\| \leq \frac{\|g_i\|}{2} \, \, \right\}, 
	\quad
	\calI^{g\nearrow} \eqdef \calI^{newt} \setminus \calI^{g\searrow}, \\
	\calI^{decr} &\eqdef  \left\{i\in \calI^{g\nearrow} \, \Big|  \; \|s_i^{trial}\| \geq \frac{1}{\sqrt{\sigma_{i}}\kap{slow}} \right\},
\end{align*} 
the last subset containing the indices of the Newton iterations where both conditions in \req{condreject} fail.
The corresponding subsets of successful iterations are then given by
\[
\begin{array}{lclclcl}
	\calS_k^{newt} &\eqdef& \calS_k \cap \calI^{newt}, 
	&&
	\calS_k^{ncurv} &\eqdef& \calS_k \cap \calI^{ncurv}\\
	\calS_k^{g\searrow} &\eqdef& \calS_k^{newt} \cap \calI^{g\searrow},
	&& \calS_k^{decr} &\eqdef& \calS_k^{newt} \cap \calI^{decr}.
\end{array}
\]
Since the iteration is unsuccessful if the test \eqref{condreject} holds, one checks that
\beqn{Skneigdivision}
\calS_k^{newt} \eqdef \calS_k^{g\searrow} \cup \calS_k^{decr}.
\eeqn 

We also recall a well-known result bounding the total number of
iterations of adaptive regularization methods in
terms of the number of successful ones.

\llem{SvsU}{
	\cite[Theorem~2.4]{BirgGardMartSantToin17},\cite[Lemma~2.4.1]{Cartis2022-wb}
	Suppose that the \al{AN2CLS} algorithm is used and that $\sigma_k \leq
	\sigma_{\max}$ for some $\sigma_{\max} >0$. Then
	\beqn{unsucc-neg}
	k \leq |\mathcal{S}_k| \left(1+\frac{|\log\gamma_1|}{\log\gamma_2}\right)+
	\frac{1}{\log\gamma_2}\log\left(\frac{\sigma_{\max}}{\sigma_0}\right).
	\eeqn
}

This result implies that the overall complexity of the algorithm can be
estimated once bounds on $\sigma_k$ and $|\mathcal{S}_k|$ are known, as we will show in the
next section.

We now state some upper bounds on the stepsize in all cases.

\llem{skbound}{Let $k$ an iteration of \tal{AN2CLS}. Then,
	\beqn{skneigbound}
	\|s_k^{trial}\| \leq  \sqrt{\frac{1}{\sigma_{k}}}
	\tim{ for} k \in \calI^{newt}
	\eeqn
	and
	\beqn{skcurvbound}
	\|s_k^{trial}\| = \frac{\theta \kappa_C}{\sqrt{\sigma_k}} \tim{ for } k \in \calI^{ncurv}.
	\eeqn
}
\proof{Let $k \in \calI^{newt}$. 
	Using the definition of $r_k^{trial}$ in \eqref{rneig-def}, the Cauchy-Schwartz and \eqref{orthcond}, we derive that
	\[
	(s_k^{trial})^\intercal ( H_k + (\sqrtsigkgk + \mu_k)I_n) s_k^{trial}  
	= -(s_k^{trial})^\intercal g_k +  (s_k^{trial})^\intercal r_k^{trial} 
	\leq \|s_k^{trial}\| \|g_k\|,
	\]
	But \eqref{posdefnewton} implies that
	\[
	\sqrt{\sigma_k} \,\|g_k\|\,\|s_k^{trial}\|^2 
	\leq (s_k^{trial})^\intercal ( H_k + (\sqrtsigkgk + \mu_k)I_n) s_k^{trial},
	\]
	and \eqref{skneigbound} thus follows.
	The second equation results from \eqref{sncurv-def} and \eqref{negcurvvector}.
}

\noindent
It results from this lemma that the  bound on local Lipschitz error stated in Lemma~\ref{Lipbounds} holds for $k \in \calI^{newt}$ with $\sigma_{k} \geq \frac{1}{\delta^2}$ or for $k \in \calI^{ncurv}$ with $\sigma_{k} \geq \frac{\theta \kappa_C}{\delta^2}$.

We now give a lower bound on the local second-order Taylor approximation at iteration $k$ depending on the step's type.

\llem{quaddecrease}{Let $k$ be an iteration of \tal{AN2CLS}. Then, 
	\beqn{quaddecrneig}
	-(g_k^\intercal s_k + \frac{1}{2} s_k^\intercal H_k s_k ) \geq  \sqrtsigkgk \|s_k\|^2
	\tim{ for} k \in \calI^{newt}
	\eeqn
	and
	\beqn{quaddecrcurv}
	-(g_k^\intercal s_k + \frac{1}{2} s_k^\intercal H_k s_k ) \geq \frac{1}{2} \theta^3 \kappa_C^3 \frac{\|g_k\|}{\sqrt{\sigma_{k}}} = \frac{1}{2} \sigma_{k}  \|g_k\| \|s_k\|^3
	\tim{ for } k \in \calI^{ncurv}.
	\eeqn
}
\proof{
	Consider first  the case where $k \in \calI^{newt}$. By using the definition of $r_k^{trial}$ in \eqref{rneig-def}, \eqref{orthcond} and  \eqref{posdefnewton}, we obtain that
	\begin{align*}
		g_k^\intercal s_k^{trial} + \frac{1}{2} (s_k^{trial})^\intercal H_k s_k^{trial} 
		&= (r_k^{trial})^\intercal s_k^{trial} - \frac{1}{2} (s_k^{trial})^\intercal H_k s_k^{trial} - \sqrtsigkgk \|s_k^{trial}\|^2 - \mu_k \|s_k^{trial}\|^2 \\
		&=  -\frac{1}{2} (s_k^{trial})^\intercal (H_k + \mu_kI_n ) s_k^{trial}
		-  \sqrtsigkgk \|s_k^{trial}\|^2 - \frac{1}{2} \mu_k \|s_k^{trial}\|^2 \\
		&\leq -\sqrtsigkgk \|s_k^{trial}\|^2,
	\end{align*}
	yielding \eqref{quaddecrneig}.
	
	Suppose now that $k \in \calI^{ncurv}$. Since \eqref{sncurv-def} and \eqref{negcurvvector} hold  and since Step~4 of \al{AN2CLS} is executed when $\mu_k >  \kappa_C \sqrtsigkgk$, we deduce that,
	\begin{align}\label{quadcurvfp}
		g_k^\intercal s_k^{trial} + \frac{1}{2} (s_k^{trial})^\intercal H_k s_k^{trail} &\leq \frac{1}{2}  \|s_k^{trial}  \|^2 u_k^\intercal H_k u_k \leq -\frac{1}{2}  \frac{\theta^3\kappa_C^2}{\sigma_{k}}  \mu_k \leq -\frac{1}{2}  \theta^3\kappa_C^3 \frac{\|g_k\|}{\sqrt{\sigma_{k}}} 
	\end{align}
	yielding the first inequality in \eqref{quaddecrcurv}.
	The second inequality is obtained by substituting the bound
	\[
	\frac{1}{2} \theta^3 \kappa_C^3 \frac{\|g_k\|}{\sqrt{\sigma_{k}}} =  \frac{1}{2} \theta^3 \kappa_C^3 \sigma_{k} \frac{\|g_k\|}{\sigma_{k} \sqrt{\sigma_{k}}}= \frac{1}{2} \sigma_{k} \|g_k\| \|s_k^{trial}\|^3 
	\]
	in \req{quadcurvfp}, where we used \eqref{skcurvbound} to deduce the last equality.
}

\numsection{Complexity Analysis of \tal{AN2CLS}}\label{complexity-s}

We start our analysis by examining under what conditions an iteration of the \tal{AN2CLS} must be successful. We have already seen in Lemma~\ref{skbound} that $\|s_k^{trial}\|$ must be less than $\delta$ for large enough $\sigma_k$, and  now investigate the conditions of a successful iteration under this assumption.  We first 
consider the case of Newton iterations and examine the conditions under which the algorithm moves to Step~5, which is necessary if iteration $k$ is to be successful. 

\llem{skneiglowbound}{
	Suppose that AS.1 and AS.3 hold and that $\|s_k^{trial}\| \leq \delta$. If  $k \in \calI^{g\nearrow}$,  then we have that 
	\beqn{gkplusonesmall}
	\|s_k^{trial}\| \geq \frac{1}{\sqrt{\sigma_{k}} \left((1+\kappa_\theta+\kappa_C) + \sqrt{(1+\kappa_\theta+\kappa_C)^2 + \frac{\frac{L_0}{\|g_k\|} + L_1}{\sigma_{k}} }\right)}.
	\eeqn
	Moreover, the algorithm proceeds to Step~5 for all $k \in \calI^{newt}$ provided
	\beqn{bigsigma_1}
	\sigma_k \ge \frac{1}{\vartheta}\left(\frac{L_0}{\|g_k\|}+L_1\right).
	\eeqn
}
\proof{
	Consider $k \in \calI^{g\nearrow}$ such that $\|s_k^{trial}\| \leq \delta$. 
	Using \eqref{Lipg}, the fact that \\ $\|\nabla_x^1f(x_k+s_k^{trial})\| {>} \frac{\|g_k\|}{2}$ since $k \in \calI^{g\nearrow}$, \eqref{newtcond} and the inequality $ \mu_k \leq \kappa_C \sqrtsigkgk$, we derive that
	\begin{align*}
		\frac{\|g_k\|}{2}&<\|\nabla_x^1f(x_k+s_k^{trial})\| \leq \|\nabla_x^1f(x_k+s_k^{trial}) - g_k- H_ks_k^{trial}\| + \|g_k+H_ks_k^{trial}\| \\ 
		&\leq \frac{(L_0 + L_1 \|g_k\|)}{2} \|s_k^{trial}\|^2 + (\sqrtsigkgk + \mu_k) \|s_k^{{trial}}\| + \|r_k^{trial}\| \\
		&\leq \frac{(L_0 + L_1 \|g_k\|)}{2} \|s_k^{trial}\|^2 + (1+\kappa_C) \sqrtsigkgk \|s_k^{{trial}}\| + \kappa_\theta \sqrtsigkgk \|s_k^{trial}\|. 
	\end{align*}
	This gives a {quadratic} inequality in $\|s_k^{trial}\|$ whose solution is given by
	\[
	\|s_k^{trial}\| \geq \frac{-(1+\kappa_\theta+\kappa_C) \sqrtsigkgk + \sqrt{(1+\kappa_\theta+\kappa_C)^2 \sigma_{k} \|g_k\|^2 + \|g_{k}\| (L_0 + L_1 \|g_k\| )}}{ (L_0 + L_1 \|g_k\|)}.
	\]
	Taking the conjugate in the last inequality and then factorizing $\sqrtsigkgk$ in the resulting denominator gives \req{gkplusonesmall}. A comparison of this bound with the second part of \req{condreject} and the definition of $\kap{slow}$ in Step~0 shows that the second part of \req{condreject} must fail if \req{bigsigma_1} holds.  Thus the algorithm proceeds to Step~5 if $k \in \calI^{g\nearrow}$ and $\|s_k^{trial}\| \leq \delta$.
	Now, if $k\in \calI_k^{g\searrow} = \calI^{newt}\setminus\calI^{g\nearrow}$, then the first part of \eqref{condreject} fails and the algorithm also proceeds to Step~5. Thus {it} does so for all $k\in \calI^{newt}$.
}

\noindent
Assuming now that Step~5 is reached, iteration $k$ can only be successful if \req{nextgradcond} fails. The condition under which this must happen is the object of the next lemma.

\llem{gradboundnext}{Suppose that AS.1 and AS.3 hold and let $k$ such that $\|s_k\|\leq \delta$. Then, 
	\beqn{gkplusoneneig}
	\|\nabla_x^1 f(x_k+s_k^{{trial}})\| \leq \left(\frac{\frac{L_0}{\|g_k\|} + L_1}{2\sigma_k} + 1+\kappa_C+ \kappa_\theta \right) \|g_k\| \tim{ for } k \in \calI^{newt}
	\eeqn
	and	
	\beqn{gkplusonecurv}
	\|\nabla_x^1 f(x_k+s_k^{{trial}})\| \leq \left(\kappa_C^2 \theta^2\frac{\frac{L_0}{\|g_k\|} + L_1}{2\sigma_k} + 1 + \frac{\kappa_C \mu_k}{\sqrt{\sigma_{k} }} \right) \frac{\|g_k\|}{\epsilon} \tim{ for } k \in \calI^{ncurv}.
	\eeqn
	Moreover, condition \req{nextgradcond} fails provided
	\beqn{sigmabig_2}
	\sigma_k \geq \frac{1}{3(1-\eta_2)}\left(\frac{L_0}{\|g_k\|}+L_1\right).
	\eeqn
}

\proof{
	First note that  Lemma~\ref{Lipbounds} holds because $\|s_k\| \leq \delta$.  
	Let $k \in \calI^{newt}$. From the triangular inequality, \eqref{Lipg}, \eqref{newtcond}, 
	the fact that $\mu_k \leq \kappa_C \sqrtsigkgk$ when $k \in \calI^{newt}$ and 
	\eqref{skneigbound}, we deduce that,
	\begin{align*}
		\|\nabla_x^1 f(x_k+s_k^{trial})\| &\leq (L_0 + L_1 \|g_k\|) \frac{\|s_k^{trial}\|^2}{2} +  \|H_k s_k^{trial} + g_k  \| \\
		&\leq (L_0 + L_1 \|g_k\|) \frac{\|s_k^{trial}\|^2}{2} +  \|(\sqrtsigkgk+\mu_k)s_k^{trial}\| + \|r_k^{trial}\|   \\ 
		&\leq (L_0 + L_1 \|g_k\|) \frac{\|s_k^{trial}\|^2}{2} + (1+\kappa_C)  \sqrtsigkgk \|s_k^{trial}\| + \kappa_\theta \|g_k\| \\ 
		&\leq \frac{L_0 + L_1 \|g_k\|}{2\sigma_k} + (1+\kappa_C + \kappa_\theta) \|g_k\|
	\end{align*}
	and \eqref{gkplusoneneig} follows.
	
	Now, if $k \in \calI^{ncurv}$, we obtain from \eqref{sncurv-def} and {\eqref{negcurvvector}} that
	\[
	\|H_ks_k^{trial} \|
	= \frac{\theta \kappa_C}{\sqrt{\sigma_k}}\|H_ku_k\|= \frac{\theta \kappa_C}{\sqrt{\sigma_k}}\sqrt{u_k^{\intercal}H_k^2 u_k}\\
	\leq  \frac{ \kappa_C \mu_k}{\sqrt{\sigma_k}}.  
	\]
	Using \eqref{Lipg}, \eqref{skcurvbound}, the last inequality and the fact that $\|g_k\| \geq \epsilon$ 
	before termination and the bound $\epsilon \leq 1$, we derive that
	\begin{align*}
		\|\nabla_x^1 f(x_k+s_k^{trial})\| &\leq (L_0 + L_1 \|g_k\|) \frac{\|s_k^{trial}\|^2}{2} +  \|H_k s_k^{trial} + g_k  \| \\ 
		&\leq (L_0 + L_1 \|g_k\|) \frac{\|s_k^{trial}\|^2}{2}{+}\frac{\kappa_C\mu_k }{\sqrt{\sigma_k}} + \|g_k\|\\
		&\leq \kappa_C^2 \theta^2 \frac{L_0 + L_1 \|g_k\|}{2\sigma_k}
		+\frac{\kappa_C \mu_k}{\sqrt{\sigma_k}} \frac{\|g_k\|}{\epsilon} + \frac{\|g_k\|}{\epsilon},
	\end{align*}
	which yields \eqref{gkplusonecurv}.
	Finally, a comparison of \eqref{gkplusoneneig} with the definition of $\kap{k}$ in \req{newtstep} and $\kap{upnewt}$ in Step~0 shows that \req{nextgradcond} fails for $k \in \calI^{newt}$ if \req{sigmabig_2} holds. Similarly, comparing {\eqref{gkplusonecurv}} with the definition of $\kap{k}$ in \req{negcurvstep} shows that \req{nextgradcond} also fails for $k \in \calI^{ncurv}$ if \req{sigmabig_2} holds.  It therefore fails for all $k \geq 0$ under this condition.
}

\noindent
We are now in a position to establish the conditions under which iteration $k$ is successful, from which we derive a crucial upper bound on the regularization parameter $\sigma_k$.

\llem{boundsigmak}{
	Suppose that AS.1 and AS.3 hold. Then, for all $k \geq 0$,
	\beqn{sigmamax}
	\sigma_{k} \leq \frac{\kap{max}}{\epsilon},
	\eeqn
	where
	\beqn{kappamax}
	\kap{max} \eqdef \gamma_3 \max \left(\sigma_0, \, \frac{\kappa_C^2 \theta^2}{\delta^2}, \, \frac{1}{\delta^2} ,\,  \frac{(L_0+L_1)}{3(1-\eta_2)}, \, \frac{L_0+L_1}{\vartheta}\right).
	\eeqn
}
\proof{
	First consider $k \in \calI^{newt}$  and suppose that
	\beqn{sigmamaxnewt}
	\sigma_{k} \geq  \max \left[ \sigma_0, \frac{1}{\delta^2}, \, \frac{(\frac{L_0}{\|g_k\|} + L_1)}{{3}(1-\eta_2)},  \frac{\frac{L_0}{\|g_k\|} + L_1}{\vartheta} \right].
	\eeqn
	so that \eqref{skneigbound} yields that $\|s_k^{trial}\| \leq \delta$ and Lemma~\ref{Lipbounds} applies.
	We may thus use Lemma~\ref{skneiglowbound} to deduce that the algorithm proceeds to Step~5.
	Now, using the Hessian tensor Lipschitz approximation error stated in \eqref{Lipf},  
	\eqref{quaddecrneig}  and \eqref{skneigbound}, we obtain that
	\begin{align}\label{rhoskn}
		1-\rho_k 
		= \frac{f(x_k + s_k) - f(x_k) -g_k^\intercal s_k - \frac{1}{2} s_k^\intercal H_k s_k }
		{-g_k^\intercal s_k - \frac{1}{2} s_k^\intercal H_k s_k} 
		&\le \frac{(L_0 + L_1 \|g_k\|) \|s_k^{trial} \|^3}{6\, \sqrtsigkgk \|s_k^{trial} \|^2} \nonumber \\ 
		&= \frac{(L_0 + L_1 \|g_k\|)  \|s_k^{trial} \|}{6\, \sqrtsigkgk } \nonumber \\  
		&\leq \frac{(\frac{L_0}{\|g_k\|} + L_1 ) }{6\, \sigma_k} .
	\end{align}
	Hence, if \req{sigmamaxnewt} holds, then $\rho_k \geq \eta_2$.  Moreover, \req{sigmamaxnewt} also ensures that 
	\req{nextgradcond} fails, as shown in Lemma~\ref{gradboundnext}. As a consequence, iteration $k$ is successful.
	
	Consider now the case where $k \in \calI^{ncurv}$, in which case Step~5 of the algorithm is always executed. Suppose that
	\beqn{sigmamaxcurv}
	\sigma_{k} \geq  \max \left[ \sigma_0, \frac{\kappa_C^2 \theta^2}{\delta^2}, \, \frac{\frac{L_0}{\|g_k\|} + L_1}{3(1-\eta_2)} \right], 	
	\eeqn
	again ensuring, because of \eqref{skcurvbound}, that $\|s_k^{trial}\| \leq \delta$ and that the Lipschitz error bound \eqref{Lipf} holds. Using this bound and \eqref{quaddecrcurv}, we derive that
	\[
	1-\rho_k         
	= \frac{f(x_k + s_k) - f(x_k) -g_k^\intercal s_k - \frac{1}{2} s_k^\intercal H_k s_k }
	{-g_k^\intercal s_k - \frac{1}{2} s_k^\intercal H_k s_k} 
	\leq \frac{(L_0 + L_1 \|g_k\|)  \|s_k^{trial} \|^3}{6 \frac{1}{2} \sigma_{k} \|g_k\| \|s_k^{trial} \|^3} 
	\leq \frac{\frac{L_0}{\|g_k\|} + L_1}{3 \sigma_{k}}.
	\]
	Thus, \req{sigmamaxcurv} ensures that $\rho_k \geq \eta_2$.  We then apply Lemma~\ref{gradboundnext} and deduce from \req{sigmamaxcurv} that \req{nextgradcond} fails and therefore that iteration $k$ is successful.
	
	As a consequence, we obtain by { taking the maximum of both \req{sigmamaxnewt} and \req{sigmamaxcurv}}, that iteration $k$ is successful and $\rho_k \ge \eta_2$ provided 
	\[
	\sigma_{k} \geq  \max \left[ \sigma_0, \frac{1}{\delta^2}, \, \frac{\kappa_C^2 \theta^2}{\delta^2}, \,\frac{(\frac{L_0}{\|g_k\|} + L_1)}{3(1-\eta_2)},  \frac{\frac{L_0}{\|g_k\|} + L_1}{\vartheta} \right].
	\]
	The update formula \req{sigmakupdate} and the facts that $\|g_k\| \geq \epsilon$ before termination and that $\epsilon \leq 1$ then imply \req{sigmamax}-\req{kappamax}.
}

We now provide also an upper-bound on $|\calS_k^{g\searrow}|$ in the spirit of \cite[Lemma~3.4]{GratJerToin23bIMA}. This is where  the test \eqref{nextgradcond} is needed to bound the ratio $\frac{\|g_{k+1}\|}{\|g_k\|}$ for all iterations, irrespective of  the Lipschitz error bounds \eqref{Lipg}. 

\llem{gdecr}{Suppose that AS.1, AS.3 and AS.4 hold. Then,
	\beqn{upperSkdg}
	|\calS_k^{g\searrow}| \leq  \log\left( \frac{\kap{upnewt}}{{\epsilon}}\right) \frac{|\calS_k^{decr}|}{\log(2)}  
	+ \log\left(\frac{\kap{upncurv}}{{\epsilon}}\right) \frac{|\calS_k^{ncurv}|}{\log(2)} +   \frac{|\log(\epsilon)| + \log(\|g_0\|)}{\log(2)} + 1,
	\eeqn
	where $\kap{upnewt}$ is defined in Step~0 and $\kap{upncurv}$ is given by
	\beqn{kapupncurvdef}
	\kap{upncurv} \eqdef \frac{3}{2} \theta^2 \kappa_C^2(1-\eta_2) + 1 + \frac{\kappa_C \kappa_B}{\sqrt{\sigma_{\min}}}.
	\eeqn
} 

\proof{
	First observe that if $k \in \calS_k^{g\searrow}$, $\|g_{k+1}\| \leq \frac{\|g_k\|}{2}$. Let $k \in \calS_{ k}^{ncurv}$. Using \eqref{nextgradcond} with $\kappa_k$ defined in \eqref{negcurvstep}, the facts that $\sigma_{k} \geq \sigma_{\min}$ and ${\mu_k}\leq \frac{\kappa_B}{\theta}$ from AS.4 and \eqref{negcurvvector},
	we obtain that
	\beqn{gkplusoncurvbound}
	\frac{\|g_{k+1}\|}{\|g_k\|} 
	\leq \frac{\frac{3}{2}\theta^2 \kappa_C^2(1-\eta_2)+1+\frac{\kappa_C \mu_k}{\sqrt{\sigma_{k}}}} {\epsilon} 
	\leq \frac{\frac{3}{2}\theta^2 \kappa_C^2(1-\eta_2)+1+\frac{\kappa_C \kappa_B}{\sqrt{\sigma_{\min}}}}{\epsilon} 
	= \frac{\kap{upncurv}}{\epsilon}.
	\eeqn
	Successively using the fact that $\calS_k = \calS_k^{decr}  \cup \, \calS_k^{g\searrow} \cup \, \calS_k^{ncurv}$, the relationship between $\|g_{k+1}\|$ and $\|g_k\|$ in the three cases  depending whether the iterate $i \in \calS_k^{decr} $ ($\kappa_k = \kap{upnewt}$ in \eqref{nextgradcond}), or $i \in \calS_k^{ncurv}$ \eqref{gkplusoncurvbound}, or $i \in \calS_k^{g\searrow}$, we derive that, for $k\ge 0$
	\begin{align*}
		\frac{\epsilon}{\|g_0\|} &\leq \frac{\|g_k\|}{\|g_0\|} = \prod_{i \in \calS_k \setminus \{k\}} \frac{\|g_{i+1}\|}{\|g_i\|} \\
		& = \prod_{i \in \calS_k^{decr} \setminus \{k\}} \frac{\|g_{i+1}\|}{\|g_i\|} \prod_{i \in \calS_k^{g\searrow} \setminus \{k\}} \frac{\|g_{i+1}\|}{\|g_i\|} \prod_{i \in \calS_k^{ncurv} \setminus \{k\}} \frac{\|g_{i+1}\|}{\|g_i\|} \\
		&\leq   \left(\frac{\kap{upnewt}}{\epsilon}\right)^{|\calS_k^{decr} \setminus \{k\}|} \times  \frac{1}{2^{|\calS_k^{g\searrow} \setminus \{k\}|}} \times \left(\frac{\kap{upncurv}} {\epsilon}\right)^{|\calS_k^{ncurv} \setminus \{k\}|}.
	\end{align*}  
	Rearranging the last inequality and using that $|\calS_k^{ncurv} \setminus \{k\}| \leq |\calS_k^{ncurv} |$ and that $|\calS_k^{decr} \setminus \{k\}| \leq |\calS_k^{decr} |$ gives that
	\[
	\frac{2^{|\calS_k^{g\searrow} \setminus \{k\}|} \epsilon}{\|g_0\|} \leq  \left(\frac{\kap{upnewt}}{\epsilon}\right)^{|\calS_k^{decr} |} \times \left(\frac{\kap{upncurv}} {\epsilon}\right)^{|\calS_k^{ncurv} |}. 
	\]
	Taking the logarithm in this inequality, using that $|\calS_k^{g\searrow} \setminus \{k\}| \geq |\calS_k^{g\searrow} | -1$ and further rearranging, finally yields \eqref{upperSkdg}.
}


Combining the previous lemmas, we are now able to state the complexity of the \al{AN2CLS} algorithm.

\lthm{complexity}{Suppose that AS.1--AS.4 hold. Then the \al{AN2CLS} algorithm requires at most
	\[
	|\calS_k| \leq \left(\kappa_\star + \frac{(\kap{decr} + \kap{ncurv}) |\log(\epsilon)| }{\log(2)}  \right) \epsilon^{-3/2} + \frac{|\log(\epsilon)| + \log(\|g_0\|)}{\log(2)} + 1.
	\]
	successful iterations and at most
	\begin{align*}
		&\left(1+\frac{|\log\gamma_1|}{\log\gamma_2}\right)\left[\left(\kappa_\star + \frac{(\kap{decr} + \kap{ncurv}) |\log(\epsilon)| }{\log(2)}  \right) \epsilon^{-3/2} + \frac{|\log(\epsilon)| + \log(\|g_0\|)}{\log(2)} + 1 \right] \\
		&+
		\frac{1}{\log\gamma_2}\left(\log\left(\frac{\kap{max}}{\sigma_0}\right) + |\log(\epsilon)| \right)
	\end{align*}
	iterations to produce a vector $x_\epsilon$ such that $\|g(x_\epsilon)\| \leq \epsilon$, where $\kappa_\star$ is defined by
	\beqn{kapstardef}
	\kappa_\star \eqdef \kap{decr} \left(1+\frac{\log\left(\kap{upnewt}\right)}{\log(2)} \right) + \kap{ncurv} \left(1+ \frac{\log\left(\kap{upncurv} \right)}{\log(2)} \right)
	\eeqn
	with both $\kap{upnewt}$ and $\kap{upncurv}$ defined in Step~0 and \eqref{kapupncurvdef}, respectively. In addition, $\kap{ncurv}$ and $\kap{decr}$ are defined by
	\beqn{kcurvanddecrdef}
	\kap{ncurv}
	\eqdef \frac{2 (f(x_0) - f_{\rm low})  \sqrt{\kap{max}}}{\eta_1 \kappa_C^3 \theta^3 }, \quad \quad
	\kap{decr} \eqdef \frac{(f(x_0) - f_{\rm low})  \kap{slow}^2\sqrt{ \kap{max} } }{\eta_1 },
	\eeqn
	and $\kap{max}$  is defined by \eqref{kappamax}.
}

\proof{	First note that we only need to prove an upper bound on $|\calS_k^{decr}| $ and $|\calS_k^{ncurv}|$  to derive a bound on $|\calS_k|$ since
	\beqn{calSkcurv}
	| \calS_k| = |\calS_k^{decr}  | + | \calS_k^{ncurv}| + |\calS_k^{g\searrow} |
	\eeqn
	and a bound on $|\calS_k^{g\searrow} |$ is given by \eqref{upperSkdg}.
	We start by proving an upper bound on $| \calS_k^{ncurv}|$.  Using AS.2, 
	\eqref{quaddecrcurv}, that $\|g_k\| \geq \epsilon$ before termination and \eqref{sigmamax}, we derive that,
	\begin{align*}
		f(x_0) - f_{\rm low}
		\geq \!\!\!\sum_{i \in \calS_k} f(x_i) - f(x_{i+1}) \geq  \sum_{i \in \calS_k^{ncurv}} f(x_i) - f(x_{i+1})
		&\geq \sum_{i \in \calS_k^{ncurv}}  \, \frac{\eta_1 \theta^3 \kappa_C^3 }{2 \sqrt{\sigma_{i}}}\,\|g_i\| \\
		&\geq\frac{\eta_1 \theta^3 \kappa_C^3 }{2 \sqrt{\kap{max} }} \epsilon^{3/2} |\calS_k^{ncurv}|,
	\end{align*}
	and hence after rearranging,
	\beqn{Skcurvbound}
	|\calS_k^{ncurv}|
	\leq \frac{2 (f(x_0) - f_{\rm low})  \sqrt{ \kap{max}}}{\eta_1 \kappa_C^3 \theta^3 } \,\epsilon^{-{3}/{2}}
	= \kap{ncurv}\,\epsilon^{-{3}/{2}},
	\eeqn
	with $\kap{ncurv}$ defined in \eqref{kcurvanddecrdef}.
	
	The reasoning is similar for $|\calS_k^{decr}|$. Using AS.2, \eqref{quaddecrneig}, that $\|s_i^{decr}\| \geq \frac{1}{\sqrt{\sigma_{i}} \kap{slow}}$ for $i \in \calS_k^{decr}$, that $\|g_i\|\geq \epsilon$ before termination, and  \eqref{sigmamax} yields that
	\begin{align*}
		f(x_0) - f_{\rm low}	&\geq \sum_{i \in  \calS_k^{decr}} f(x_i) - f(x_{i+1})
		\geq \sum_{i \in  \calS_k^{decr}}\eta_1  \, \sqrt{\sigma_i} \|g_i\| \|s_i\|^2 \\ &\geq \sum_{i \in \calS_k^{decr}}\eta_1     \frac{\|g_i\| }{\sqrt{\sigma_{i}} \kap{slow}^2} 
		\geq |\calS_k^{decr}|\eta_1   \frac{\epsilon^\sfrac{3}{2}}{\sqrt{\kap{max}} \kap{slow}^2}. 
	\end{align*}
	Rearranging the last inequality to upper-bound $|\calS_k^{decr}|$, 
	\beqn{Skdecrbound}
	|\calS_k^{decr}| \leq \frac{(f(x_0) - f_{\rm low})   \kap{slow}^2 \sqrt{ \kap{max} }}{\eta_1} \epsilon^{-3/2}  \\
	\leq \kap{decr} \epsilon^{-3/2},
	\eeqn
	with $\kap{decr}$ defined in \eqref{kcurvanddecrdef}.
	Combining now \eqref{Skcurvbound} and \eqref{Skdecrbound} with the upper bound in \eqref{upperSkdg} and using that $|\calS_k| = |\calS_k^{decr}| + |\calS_k^{ncurv}| + |\calS_k^{g\searrow}|$, we derive that 
	\[
	|\calS_k| \leq \left(\kappa_\star + (\kap{ncurv} + \kap{decr}) \frac{|\log(\epsilon)|}{\log(2)}\right) \epsilon^{-3/2} + \frac{|\log(\epsilon)| + \log(\|g_0\|)}{\log(2)} + 1
	\]
	where $\kappa_\star$ defined in \eqref{kapstardef}. The second statement of the Theorem is proved by using both  \eqref{sigmamax} and \eqref{unsucc-neg}.
}

\noindent
Thus the complexity of the \al{AN2CLS} algorithm is, in order, the same as that of the fast Newton method proposed in \cite{GratJerToin23bIMA}, albeit under a weaker local Lipschitz smoothness condition.

As mentioned in the introduction, this new class includes a broader class of functions, such as univariate polynomials, and thus can tackle a large class of problems compared to the standard second-order methods \cite{Cartis2022-wb}.  Note that our bound differs by a $|\log(\epsilon)|$ factor from that for optimal second-order methods when searching for an approximate first-order stationary point \cite{Carmon2019a}. 

\numsection{Finding Second-Order Critical Points}\label{so-teaser}

Can the \al{AN2CLS}  algorithm be strengthened to ensure it will compute
second-order critical points? We show in this section that this is possible under the same
assumptions as those used for the first-order analysis. 
The resulting modified algorithm, which we call \al{SOAN2CLS} (for 
\al{Second-Order AN2CLS}) makes extensive use of \al{AN2CLS}, and is
detailed \vpageref{SOAN2CLS}.

\algo{SOAN2CLS}{Second-Order Adaptive Newton with Negative Curvature Local Smoothness  (\tal{SOAN2CLS})}{
	\begin{description}
		\item[Step 0: Initialization: ] 
		Identical to $\al {AN2CLS}[{\rm Step~0}]$ with $\epsilon\in(0,1]$ now replaced 
		by $\epsilon_1 \in (0,1]$ and $\epsilon_2 \in (0,1]$.
		\item[Step 1: Compute current derivatives: ]
		If not available, evaluate $g_k$ and $H_k$. Terminate if
		\beqn{SO-termination}
		\|g_k \| \leq \epsilon_1 \tim{and} \lambda_{\min}(H_k) \geq - \epsilon_2.
		\eeqn
		\item[Step 2: Step calculation: ]
		If $\| g_k \| > \epsilon_1$,
		compute a first-order step $s_k \eqdef s_k^{fo}$, $\kappa_k$ and the boolean variable \texttt{REJECT} as in Steps~2 to 4 of Algorithm~\ref{AN2Cgen} using the \tal{Stepcomp} procedure.
		Otherwise (i.e. if $\|g_k \| \leq \epsilon_1$), compute $u_k$ such that 
		\beqn{negcurvsecondorder}
		g_k^\intercal u_k  \leq 0, \, \,  \|u_k\| = 1
		\tim{and}
		H_k u_k = \lambda_{\min}(H_k) u_k,  \, \, 
		\eeqn
		and set 
		\beqn{negcurvstephess}
		s_k = s_k^{so} \eqdef  \frac{ u_k}{\sqrt{\sigma_k}}, \tim{ and } \kappa_{k,hess} \eqdef \frac{3(1-\eta_2) |\lambda_{\min}(H_k)|}{2 \sqrt{\sigma_{\min}}} + 1 + \frac{|\lambda_{\min}(H_k)|}{\sqrt{\sigma_{k}}}.
		\eeqn
		\item[Step 3: Acceptance ratio computation: ] 
		If $s_k = s_k^{fo}$, use Step~5 of the \al{AN2CLS} algorithm with $\epsilon= \epsilon_1$.
		Otherwise, evaluate $f(x_k +s_k)$ and compute the acceptance ratio $\rho_k$ as in \eqref{rhokdef}.
		If $\rho_k < \eta_1$ and
		\beqn{nextgradcondhess}
		\|\nabla_x^1 f(x_k + s_k)\| > \kappa_{k,hess}, 
		\eeqn
		set \texttt{REJECT = TRUE}. 
		\item[Step~4: Variables' update: ]		
		If \texttt{REJECT = FALSE}, set $x_{k+1} = x_k + s_k$, otherwise set $x_{k+1} = x_k$.
		\item[Step 5: Regularization parameter update: ] 
		Identical to Step~7 of the \al{AN2CLS} algorithm.
	\end{description}
}

Before reaching an approximate
first-order point, the algorithm only uses the \tal{Stepcomp}
subroutine to generate trial steps, hence the '$fo$' (first-order)
superscripts in the first condition of Step~2. Once an approximate
first-order point is reached, further progress towards second-order
criticality is obtained by exploiting the negative-curvature
direction \eqref{negcurvsecondorder}-\eqref{negcurvstephess}, which
justifies the '$so$' (second-order) superscript.
The test \eqref{nextgradcondhess} is required to ensure that the gradient remains bounded when switching from a second-order step to a first-order one. This will be crucial in order to derive an equivalent of Lemma~\ref{gdecr} for the \al{SOAN2CLS} Algorithm.

An upper bound on the evaluation complexity of the \al{SOAN2CLS}
algorithm is given by the following theorem.

\lthm{ComplexityboundSecond}{Suppose that AS.1--AS.4 hold. Then the \al{SOAN2CLS} algorithm requires at most
	\[
	|\calS_{ k} | \leq \kappa_{\star,hess} \left(\epsilon_1^{-3/2} + \epsilon_2^{-3}\right)  + \kappa_{\star,log} |\log(\epsilon_1)| (\epsilon_1^{-3/2} + \epsilon_2^{-3}) + \frac{|\log(\epsilon_1)| + \log(\kap{gpi})}{\log(2)} + 1
	\]
	successful iterations and at most
	{ \small
		\begin{align*}
			&\left(1 + \frac{|\log(\gamma_1)|}{\log(\gamma_2)}\right) \left(\kappa_{\star,hess} \left(\epsilon_1^{-3/2} + \epsilon_2^{-3}\right)  + \kappa_{\star,log} |\log(\epsilon_1)| (\epsilon_1^{-3/2} + \epsilon_2^{-3}) + \frac{|\log(\epsilon_1)| + \log(\kap{gpi})}{\log(2)} + 1 \right) \\
			&+ \frac{1}{\log \gamma_2} \left( \log\left( \max\left(\frac{\kap{max1} \epsilon_1^{-1}}{\sigma_0}, \, \frac{\kap{max2}\epsilon_2^{-2}}{\sigma_0}\right) \right)\right)
		\end{align*}
	}
	iterations to produce a vector $x_\epsilon$ such that $\|g(x_\epsilon)\| \leq \epsilon_1$ and $\lambda_{\min}(H_k) \geq -\epsilon_2$ where $\kappa_{\star,hess}$ is defined by 
	\beqn{kappastarhess}
	\kappa_{\star,hess} \eqdef \kap{ncurvhess} \left(1+ \frac{\log\left(\kap{upncurv}\right) }{\log(2)} \right)
	+ \kap{decrhess} \left(1 + \frac{\log(\kap{upnewt})}{\log(2)}\right) + \kap{so} \left(2+\frac{\log(\kap{gpi})}{\log(2)}\right),
	\eeqn
	$\kappa_{\star,log}$ and $\kap{so}$ by
	\beqn{kappastarlogkapso}
	\kappa_{\star,log} \eqdef \frac{\kap{ncurvhess} + \kap{decrhess} + \kap{so}}{\log(2)}, \ms \kap{so } \eqdef\frac{ 2 (f(x_0) - f_{\rm low})  \max\left( {\kap{max2}}, \kap{max1}  \right)}{\eta_1 }, 
	\eeqn 
	$\kap{ncurvhess}$ and $\kap{decrhess}$ by
	\beqn{kaphessdef}
	\kap{ncurvhess} \eqdef  \kap{ncurv} \max(1, \, \frac{\sqrt{\kap{max2}}}{\sqrt{\kap{max1}}}), \ms  \kap{decrhess} \eqdef \kap{decr} \max(1, \, \frac{\sqrt{\kap{max2}}}{\sqrt{\kap{max1}}}),
	\eeqn
	and
	$\kap{max1} \eqdef \kap{max}$ is defined in \eqref{kappamax}, $\kap{max2}$ in \eqref{kappamax2}, and $\kap{gpi}$ in \eqref{gpibound}.
}

Note that the complexity bound stated in the theorem differs from that of  standard second-order optimal methods \cite[Theorems~3.3.9 and 3.4.6]{Cartis2022-wb} by a logarithmic factor $|\log(\epsilon_1)|$. However, this result is obtained using only a local smoothness assumption, while the above references typically require stronger global Lipschitz smoothness.
To prove Theorem~\ref{ComplexityboundSecond}, several modifications of the first-order theory must be considered. First, the bound of Lemma~\ref{boundsigmak} is no longer valid since we also consider a second-order step. Since the bound of $\sigma_{k}$ for \al{AN2CLS} \eqref{sigmamax} depends on $\epsilon_1$, the new bound of $\sigma_{k}$ for \al{AN2CLS} will depend on both $\epsilon_1$ and $\epsilon_2$ (see \eqref{sigmamaxhess}). Also note that Lemma~\ref{gdecr} is no longer valid, since an iterate $x_i$ with $\|g(x_i)\| \leq \epsilon_1$ can occur.  However, the required modifications remain in the spirit of the proof in Section~\ref{complexity-s}, and details of the proof of Theorem~\ref{ComplexityboundSecond} have been moved to Appendix~{\ref{app:so-theory}.}

\numsection{Step Computation and Preconditioning}\label{stepcomproutine}

In this section, we propose two methods for computing the trail step $s_k^{trial}$ in Step~2 of the \al{AN2CLS} algorithm that also satisfy the requirements of {AS.0}.  We also discuss how a preconditioned variant can be obtained. A first variant exploits exact negative curvature information and exact solution of the Newton system  to compute the step {and} a second uses Krylov subspaces of increasing dimension. Because of our focus on a particular iteration $k$, we will drop the subscript for the remainder of the section. 

\subsection{Exact Step Computation}

The first considered variant relies on exact information {and} yields an algorithm   named \al{AN2CLSE}. We now detail its step computation procedure in the \al{StepcompExact} algorithm.

\algo{StepCompExact}{ $[ s^{trial}, \mu ]$ = \al{StepcompExact}$(\, g,\, H, \, \sigma,\, \kappa_C, \, {0}, \, 1)$}{
	Set $\mu = [-\lambda_{\min}(H)]_+$.
	If $\mu \leq \kappa_C \sqrt{\sigma}\|g\| $, then return $[s^{trial},\mu]$, where $s^{trial}$ is the exact solution of the system
	\beqn{sneigexac}
	(H + (\mu + \sqrt{\sigma} \|g\|)I_n ) s^{trial} = -g.
	\eeqn
	Otherwise (i.e.\ if $\mu  > \kappa_C \sqrt{\sigma}\|g\|$), compute the eigenvector $u$ satisfying
	\beqn{uvecexact}
	g^\intercal u \leq 0, \, \, \|u\| = 1, \, \, Hu = \lambda_{\min}(H)u
	\eeqn
	and return $[s^{trial}, \mu] = \left[\bigfrac{\kappa_C}{\sqrt{\sigma}}u,\mu\right]$.
}

\noindent
We now show that this procedure qualifies for computing the trial step in \al{AN2CLS} algorithm, yielding the 
\al{AN2CLSE} variant.

\llem{exacthm}{The output of the \al{StepcompExact} procedure satisfies Assumption~\ref{assump0} with $\theta = 1$ and $\kappa_\theta=0$.}
\proof{Consider first the case where $\mu = -\lambda_{\min}(H) \le \kappa_C \sqrt{\sigma}\|g\|$. Since $\lambda_{\min}(H + \mu I_n) \geq 0$, \eqref{posdefnewton} holds. And since $s^{trial}$ is the exact solution of \eqref{sneigexac}, \eqref{newtcond} holds with $\kappa_\theta=0$ and \eqref{orthcond} also holds because $r_k^{{trial}} = 0$. 
	Else if $\mu > \kappa_C \sqrt{\sigma}\|g\|$, the exact eigenvector $u$ in \eqref{uvecexact} obviously satisfies \eqref{negcurvvector} with $\theta=1$.
}

\subsection{Krylov Variant}

When the dimension of the problem grows and factorizations, the workhorse of exact solution of linear systems, become impractical, one can turn to exploiting Krylov subspaces, as we now
show. The resulting algorithmic variant
will be called \al{AN2CLSK}, where \al{K} stands for Krylov, and is obtained by replacing Step~2 of
the \al{AN2CLS} algorithm by {Algorithm~\al{StepcompK}} \vpageref{AN2CKStep}.

\algo{AN2CKStep}{$[s^{trial},\mu]$ = \al{StepcompK}$(\,g, H, \sigma, \kappa_C,  \kappa_\theta, \theta \,)$ }{
	\begin{description}
		\item[Step 0: Initialization: ] Set $p=1$, $r_1 = g$, $\alpha_1=\|g\|$ and $z_{0} = 0$.
		\item[Step 1: Form the orthonormal basis: ] Compute
		\beqn{lanczoquantities1}
		v_p = \frac{r_p}{\alpha_p}, \ms
		\delta_p = v_p^\intercal H v_p, \ms
		\eeqn
		\beqn{lanczoquantities2}
		r_{p+1} = H v_p - \delta_p v_p - \alpha_p v_{p-1},\ms
		\alpha_{p+1} = \|r_{p+1}\|,
		\eeqn
		and define 
		\beqn{Vpdef}
		V_p = (v_1 , v_2, \dots , v_p) \in \Re^{n \times p}.
		\eeqn
		
		\item[Step 2: Newton step computation: ] 
		Form the subspace Hessian
		\beqn{Tpcompute}
		T_p \eqdef V_p^\intercal H V_p = \begin{pmatrix}
			\delta_1 & \alpha_2  &    &      &  \\
			\alpha_2       & \delta_2  & \alpha_3 &       &    \\
			&   \ddots      & \ddots         & \ddots    &     \\
			&         &           &     \delta_{p-1}      & \alpha_p   \\
			&                 &  & \alpha_{p}           & \delta_p
		\end{pmatrix},
		\eeqn
		compute its minimum eigenvalue and set $\mu = \max(0,-\lambda_{\min}(T_p))$.  
		If $\mu >  \kappa_C \sqrt{\sigma} \|g \|$, go to Step 4. 
		Otherwise, solve 
		\begin{equation}\label{steplargegradlanczos}
			\left( T_p + (\sqrt{\sigma} \| g\| +  \mu ) I_p \right) y_p = - \alpha_1 e_1.
		\end{equation}
		
		\item[Step 3: Check the linear residual: ] If
		\beqn{resdidtestsubspace}
		|\alpha_{p+1} (e_p^{\intercal} y_p)| \leq \kappa_\theta \min( \sqrt{\sigma}\|g\| \| y_p\|, \,  \|g\| ),
		\eeqn
		then return
		\beqn{skneiglanczos}
		[s^{trial},\mu]  = [ V_p y_p, \mu ].
		\eeqn
		Else increment $p$ by one and go back to Step~1.
		
		\item[Step 4: Negative curvature step: ] 
		Compute $u_p$ such that 
		\beqn{negcurvvectorlanczos}
		e_1^\intercal u_p  \leq 0, \, \,  \|u_p\| = 1, \, \,
		u_p^\intercal T_p u_p \leq \theta \,\lambda_{\min}(T_p) 	\tim{and}
		{u_p^\intercal T_p^2 u_p \leq \frac{\lambda_{\min}(T_p)^2}{2\theta^2}}.
		\eeqn
		\item[Step 5: Check quality of eigenvector] If
		\beqn{residnegcurvstep}
		|\alpha_{p+1} (e_p^\intercal u_p) |^2 \leq \frac{\lambda_{\min}(T_p)^2}{2 \theta^2} 
		\eeqn
		then return 
		\beqn{negcurvsteplanczos}
		[s^{trial},\mu] =  \left[\frac{\theta \kappa_C}{\sqrt{\sigma}}\, V_p u_p, \mu \right].
		\eeqn
		Else increment $p$ by one and go back to Step~1.
	\end{description}
}

{First, observe that since $\theta \leq 1$, $u_p$ can be chosen as the negative curvature of $T_p$. We introduced the conditions in \eqref{negcurvvectorlanczos} to allow more flexibility in the choice of $u_p$ and exploit previous Krylov subspace information. }

Each iteration of the {\al{StepcompK}} algorithm has a moderate cost
(a few vector assignments, one matrix-vector product, and --possibly-- the computation of the
smallest eigenvalue of a tridiagonal matrix {when the computation of $u_p$ is required}, see \cite{CoakRokh13} and
the references therein for details).  We observe that
\req{lanczoquantities1}-\req{lanczoquantities2} amounts to using the
standard  Lanczos process
for building an orthonormal basis $V_p$ of successive Krylov
subspaces. For a more detailed discussions on Krylov methods applied for the regularized Newton method satisfying Assumption~\ref{assump0}, see \cite[Subsection~5.2]{GratJerToin23bIMA}.
Observe that the algorithm always terminates since when $p = n$, either a negative curvature step associated with $\lambda_{\min}(H)$ is computed and  satisfies both \eqref{negcurvvectorlanczos} and \eqref{residnegcurvstep} since $\theta \leq 1$ and $\alpha_{n+1} = 0$ or $\alpha_{n+1} = 0$ and \eqref{resdidtestsubspace} holds.

We now verify that Algorithm~\al{AN2CLSK} is a valid
instantiation of Algorithm~\al{AN2CLS}.

\llem{EqLanczosgeneric}{
	The output of the \al{StepcompK} algorithm satisfies Assumption~0.
}

\proof{
	We deduce from \req{lanczoquantities1}, \req{lanczoquantities2} and {established Krylov Lanczos property, see \cite[Algorithm~5.2.1]{Conn2000}}   that
	\beqn{VpHTpV}
	H V_p
	= V_p T_p + \alpha_{p+1} v_{p+1} e_p^\intercal
	= V_p T_p + \alpha_{p+1}  v_{p+1} e_p^\intercal.
	\eeqn
	{Since $( v_1 , v_2, \dots , v_p)$ is an orthonormal basis (again see \cite[Subsection~5.2.1]{Conn2000}), then $V_p^{\intercal}  v_{p+1} = 0$ which yields \eqref{Tpcompute}.}
	Note also that as $v_1 = \frac{r_1}{\|g\|}$ from \eqref{lanczoquantities1} and 
	$V_p^\intercal V_p = I_p$, 
	\beqn{V1prop}
	V_p^\intercal g = \alpha_1 V_p^\intercal v_1 = \alpha_1 e_1.
	\eeqn
	This last identity with the fact that $T_p = V_p^\intercal H V_p$ ensures that
	\eqref{steplargegradlanczos} combined with  \eqref{skneiglanczos} is a  reformulation of
	\eqref{posdefnewton} since $\mu = \max[0,-\lambda_{\min}(T_p)]$. By {using the same arguments as before  }, the first three identities of \eqref{negcurvvectorlanczos} are also equivalent with the first three of \eqref{negcurvvector} {by respectively using \eqref{V1prop}, that $V_p$ is an orthonormal basis and the identity $T_p = V_p^\intercal H V_p$.} 
	We now prove that combining the last inequality of \eqref{negcurvvectorlanczos} with \eqref{residnegcurvstep} yields the last property of \eqref{negcurvvector}. Multiplying the matrix \eqref{VpHTpV} by {$u_p$ to the right and $u_p^\intercal V_p$ to the left}, taking the squared norm and using that $V_p^\intercal  V_p = I_p$ and $V_p^\intercal  v_{p+1} = 0$ yields that 
	\begin{align*}
		u_p^\intercal V^\intercal H^2 V_p u_p
		= \|H V_p u_p\|^2 &= \|V_p T_p u_p\|^2 + \alpha_{p+1}^2 (e_p^\intercal u_p)^2 \\
		&= u_p^\intercal T_p^2 u_p + \alpha_{p+1}^2 (e_p^\intercal u_p)^2
		\leq \frac{\lambda_{\min}(T_p)^2}{\theta^2}
	\end{align*}
	where we used both \eqref{negcurvvectorlanczos} and \eqref{residnegcurvstep} to obtain the last inequality.
	Hence \eqref{negcurvvector} holds for $V_pu_p$ computed in \al{StepcompK}.
	
	We now prove that \eqref{resdidtestsubspace} implies \eqref{newtcond}.
	Using \eqref{steplargegradlanczos}, \eqref{V1prop}, \eqref{VpHTpV}, we obtain that
	\begin{align}
		H s + g
		&= H V_p y_p + \alpha_1  v_1  = H V_p y_p + \alpha_1   V_p e_1 \nonumber \\
		&= H V_p y_p -  V_p T_p y_p - (\sqrt{\sigma} \|g \| + [-\lambda_{\min}(T_p)]_+ )  V_p y_p \nonumber \\
		&= \alpha_{p+1} (e_p^{\intercal} y_p)  v_{p+1} - (\sqrt{\sigma} \|g \| + \mu )  V_p y_p \label{residexpr}.
	\end{align}
	{Rearranging the last inequality, using that $s = V_p y_p$ and taking the norm yields that}
	\begin{align*}
		\|Hs+g + (\sqrt{\sigma} \|g \| + \mu ) s\|
		&= \|Hs+g + (\sqrt{\sigma} \|g \| + \mu ) V_p y_p\| \\
		&= \|\alpha_{p+1}(e_p^\intercal y_p)v_{p+1}\| 
		= |\alpha_{p+1}(e_p^\intercal y_p)|,
	\end{align*}
	and since $\|s^{trial}\| = \|V_p y_p\| = \|y_p\|$, \eqref{resdidtestsubspace} implies \eqref{newtcond}.
	{Rearranging} now \eqref{residexpr} to express $r^{trial}$ and using that $V_p^\intercal v_{p+1} = 0$ with $s = V_p y_p$ yields \eqref{orthcond}. 
}

\subsection{Preconditioned Variants}

A preconditioner is often used to reduce the number of iterations of the Krylov methods. Even though  not explicitly included in Algorithm~\ref{AN2Cgen}, it can nevertheless be incorporated in the algorithm by considering $\sqrt{x^\intercal Mx}$ as the primal norm and $\sqrt{x^\intercal M^{-1}x}$ as the associated dual norm when solving the minimization problem \eqref{minf}. The newly required Assumption~\ref{assump0} would therefore be written as follows.

\setcounter{assumption}{-1}
\assump{assump0precond}{
	
	\noindent
	The \tal{Stepcomp} subroutine computes a tentative regularization parameter $\mu_k$ and a trial step $s_k^{trial}$ {satisfying the} following conditions.\\
	If $\mu_k \leq \kappa_C \sqrt{\sigma_k} \|g_k\|_{M^{-1}}$, then the computed $s_k^{trial}$ must be such that
	\begin{align}
		(s_k^{trial})^\intercal (H_k + \mu_k M) s_k^{trial} &\geq 0 \label{posdefnewtonMcond} \\
		\|r_k^{trial}\|_{M^{-1}} &= \|(H_k + (  \sqrt{\sigma_k} \|g_k\|_{M^{-1}}+ \mu_k)M)s_k^{trial} + g_k\|_{M^{-1}} \nonumber \\ 
		&\leq \kappa_\theta \min\left( \sqrt{\sigma_k} \|g_k\|_{M^{-1}} \|s_k^{trial}\|_{M},  \|g_k\|_{M^{-1}}\right), \label{newtcondMcond} \\ 
		(r_k^{trial})^\intercal s_k^{trial} &= 0. \label{orthMcond}
	\end{align}
	Else if $\mu_k > \kappa_C \sqrt{\sigma_k} \|g_k\|_{M^{-1}}$, $s_k^{trial}$ is given by
	$s_k^{trial} = \bigfrac{\theta \kappa_C}{\sqrt{\sigma_k}} u_k$ where the vector $u_k$  verifies
	\beqn{negcurvvectorMcond}
	g_k^\intercal u_k  \leq 0, \, \,  \|u_k \|_M = 1, \, \,
	u_k^\intercal H_k u_k \leq -\theta \mu_k
	\tim{and}
	u_k^\intercal H_k^2 u_k \leq \frac{\mu_k^2}{\theta^2}.
	\eeqn
}

Both Algorithms~\ref{StepCompExact} and \ref{AN2CKStep} can be easily modified accordingly. For Algorithm~\ref{StepCompExact}, $\mu$ becomes $[-\lambda_{\min}(M^\sfrac{-1}{2}HM^\sfrac{-1}{2})]_+
$, while the preconditioned Lanczos method should be used for the Krylov variant, see \cite[Subsection~5.2]{Conn2000} for more details. The proof that the preconditioned Algorithms~\ref{AN2CKStep} and \ref{StepCompExact} satisfy the requirements of Assumption~\ref{assump0precond} follows the lines of Lemma~\ref{exacthm} and Lemma~\ref{EqLanczosgeneric}, respectively.  The complexity analysis of Section~\ref{complexity-s} and Section~\ref{so-teaser}  also follows when using the newly defined primal and dual norms. 

\numsection{Numerical Illustration}\label{sec:numerics}

In the next section, we provide a numerical comparison between the newly proposed \al{AN2CLS} and the algorithms developed in \al{AN2C} \cite{GratJerToin23bIMA}.  Our implementation of \al{AN2C}, does not follow the subspace implementation suggested in \cite{GratJerToin23bIMA}, but uses a variant closer to Algorithm~\ref{AN2Cgen} where conditions \eqref{nextgradcond} and \eqref{condreject} are suppressed and the appropriate gradient regularization's used 
and negative curvature is employed based on another condition\footnote{ More specifically, the subspace implementation for \al{AN2C} when computing a trial step, is avoided using a condition similar to Assumption~\ref{assump0} where Newton step is taken when $\mu_k \leq \kappa_C \sqrt{\sigma_k \|g_k\| }$ and negative curvature employed otherwise. The two conditions \eqref{orthcond} and \eqref{posdefnewton} are still imposed, while the first term of the $\min$ in \eqref{newtcond} is changed to $\sqrt{\sigma_k \|g_k\|} \|s_k\|$ for this new \al{AN2C} variant. Conditions on the unitary negative curvature  direction $u_k$ \eqref{negcurvvector} are remain unchanged but the step in \eqref{sncurv-def} is scaled by $\theta \kappa_C \sqrt{\sigma_k \|g_k\|}/\sigma_k$ instead.}. In addition, we also consider two other "baseline" methods.  The first is based on cubic adaptive regularization \al{AR2} \cite{BirgGardMartSantToin17} and the second is a trust-region method \al{TR2} \cite{Conn2000}. Both these methods use the gradient and the Hessian to construct a quadratic approximation, making them suitable for a comparison with \al{AN2C}-like methods. 

\noindent
The following set of hyperparameters is used for both variants of \al{AN2CLS}
\begin{equation}\label{AN2CLShypparam}
	\kappa_C = 10^3,  \, \,   \vartheta = 10^4, \, \gamma_1 = 0.5, \, \gamma_2 = \gamma_3 = 10, \, \eta_1 = 10^{-4}, \, \eta_2 =0.95, \, \sigma_{\min} = 10^{-8}\, \quad \sigma_0 = \frac{1}{\|g_0\|}.
\end{equation}

We illustrate the performance of our algorithm on three sets of test problems from
the freely available S2MPJ collection of CUTEst problems \cite{S2MPJ}. 

The first set contains $95$ small-dimensional problems, with dimensions ranging from $2$ to $49$,
the second contains 47 medium-dimensional problems with dimensions ranging from $50$ to $997$, while the third contains $33$ "large" problems with dimensions ranging from $1000$ to $5000$. For our numerical comparison, we will use our new implementation of the \tal{AN2C} method.

\subsection{Using Exact Linear Solves and Eigenvalue Decomposition}

In this subsection, we focus on the algorithm \al{AN2CLSE} that uses Algorithm~\ref{StepCompExact} to compute the associated $(s_k, \, u_k, \mu_k)$. In this case, we impose $\theta = 1$. As a baseline, we will use three different methods. The first is  the algorithm \cite[\al{AN2CE}]{GratJerToin23bIMA}. For its hyper-parameters, we keep those of \eqref{AN2CLShypparam} and change only $\sigma_0 = 1$. For the two other baselines, we  use an adaptive cubic regularization method \al{AR2F} and a trust-region method \al{TR2E}. For both of these methods, we solve the subproblem exactly. In \al{AR2F}, we use an exact subsolver based on the the secular equation and matrix factorization (see \cite[Chapter~9]{Cartis2022-wb} for further clarification). For the trust-region subproblem, we use the exact subsolver\footnote{https://github.com/oxfordcontrol/TRS.jl/tree/master} based on the work of \cite{Adachietal17}.  The hyperparameters of \al{AR2F} are the same as in \eqref{AN2CLShypparam} when updating the regularization parameter. For the trust-region method, the radius is decreased by a factor $\sqrt{10}$ and expanded by a factor 2 as done in \cite{GratJerToin23bIMA}. All experiments were run in Julia on a machine with AMD Ryzen 7 5000 at 3.8GHz.

We stop the algorithm when the gradient norm falls below $\epsilon = 10^{-6}$. We also stop a run when either the total number of iterations exceeds 5000 or when the cpu-time exceeds one hour for a specific instance.  

Our results are summarized using the standard performance profile \cite{DolaMoreMuns06} and two additional metrics.  Efficiency is measured, in accordance with the complexity theory, in the number of iterations (or, equivalently,
function and possibly derivatives’ evaluations): the fewer the more efficient the algorithm. We also add an additional global metric  $\pi_{\tt algo}$ following \cite{PorcToin19} which denotes $1/10$ of the
curve corresponding to $\tal{algo}$ in the performance profile, for abscissas in the interval $[1, 10]$. We also add a reliability metric $\rho_{\tt algo}$. It notes the percentage of successful runs taken on all problems in each of the three classes.

\begin{figure}\label{thbp}
	\centerline{
		\includegraphics[width= 0.5 \linewidth]{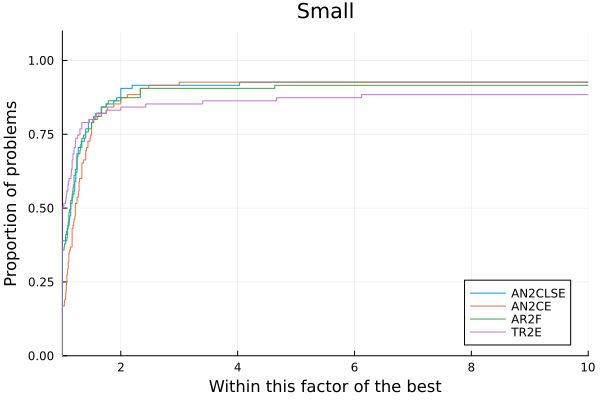}
		\hfil
		\includegraphics[width= 0.5 \linewidth]{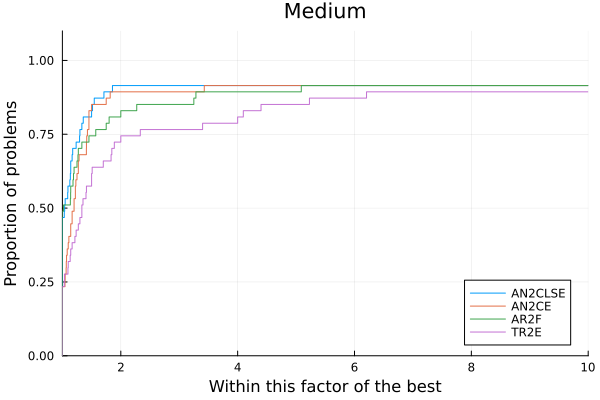}
	}
	\centerline{
		\includegraphics[width= 0.5 \linewidth]{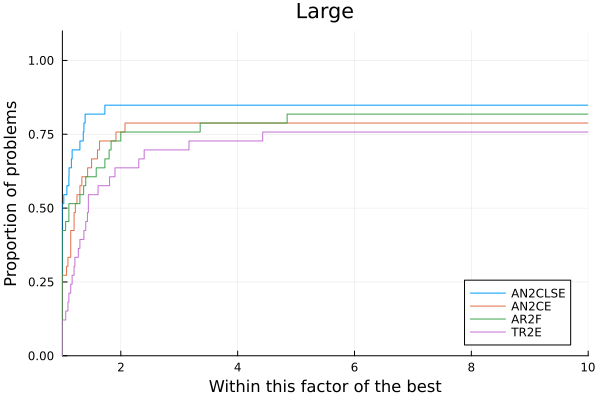}
	}
	\caption{Performance profile of both \al{AN2CLSE} and \al{AN2CE} on three different set of problems (small, medium, large)}
\end{figure}

\begin{table}[htb]\footnotesize
	\centering
	\begin{tabular}{|l|r|r|r|r|r|r|}
		\hline
		& \multicolumn{2}{|c|}{small pbs.}& \multicolumn{2}{|c|}{medium pbs.}& \multicolumn{2}{|c|}{largish pbs.}\\
		\hline
		{\tt algo} & $\pi_{\tt algo}$  & $\rho_{\tt algo}$ &
		$\pi_{\tt algo}$ & $\rho_{\tt algo}$ &$\pi_{\tt algo}$ & $\rho_{\tt algo}$ \\
		\hline
		\al{AN2CLSE}  &  0.91 &  93.68  & 0.90  &  91.48 & 0.83& 84.84 \\
		\al{AN2CE}   &  0.90 &  93.68  & 0.89  &  92.61 & 0.76 & 78.78 \\
		\al{AR2F}   &  0.89 &  91.57  & 0.88  &  91.48 & 0.79 & 81.381 \\
		\al{TR2E}   &  0.88 &  90.57  & 0.85  &  89.36 & 0.74 & 78.78 \\
		\hline
	\end{tabular}
	\caption{\label{tab:statsE} Efficiency and reliability statistics for
		the {\sf{CUTEST}} Julia problems for full space variants.}
\end{table}

We see from the results that the local smooth variant \al{AN2CLSE} performs slightly better than the other three baselines for the three different sets of problems. This might reflect the fact that \al{AN2CLSE} is designed to handle functions that satisfy a weaker local smoothness assumption  \eqref{LipHessian} than what is theoretically required for the other methods, but this clearly requires further analysis.

Note that for both variants, negative curvature steps are rarely performed (less than 1\%) despite being required theoretically making both methods  regularized Newton ones numerically (this behavior was also observed in \cite{GratJerToin23bIMA} for the \al{AN2CE} method).  

\subsection{Krylov Variants}

We now turn to variants using a Krylov iterative method in order to compute the step. We will choose additional hyperparameters for both \al{AN2CK} and \al{AN2CLSK} as 
\[
{\kappa_\theta} = 1, \quad \theta = \frac{1}{2},
\]
the other hyper-parameters being chosen as for the full space variants. The choice of ${\kappa_\theta}$ and $\theta$ was done after a search on a subset of small dimensional problems.  The choice of $\theta = \frac{1}{2}$ is to allow $u$ to be chosen as the sum of the current vector $y_p$ plus a multiple of the eigenvector associated with $\lambda_{\min}(T_p)$ chosen to ensure that the two last inequalities of \eqref{negcurvvectorlanczos} hold. This strategy was shown to be efficient in \cite{GratJerToin23bIMA}.   None of these methods uses preconditioning and the matrices $V_p$ are stored explicitly. We compare these two variants with \al{AR2K} and \al{TR2K}. For \al{AR2K}, we exactly minimize the local cubic model in an increasing Krylov subspace until $\|g_k  + H_k s_k\| \leq \sfrac{1}{2} \theta_{sub} \sigma_k \|s_k\|^2$ with $\theta_{sub}=2$, as proposed in \cite{GrattonToint22}. For the \al{TR2K} algorithm, we iteratively increase the subspace and exactly solve the subproblem until
\[
\|g_k + H_k s_k\| \leq \frac{\|g_k\|}{10}.
\]

\begin{figure}
	\centerline{
		\includegraphics[width= 0.5 \linewidth]{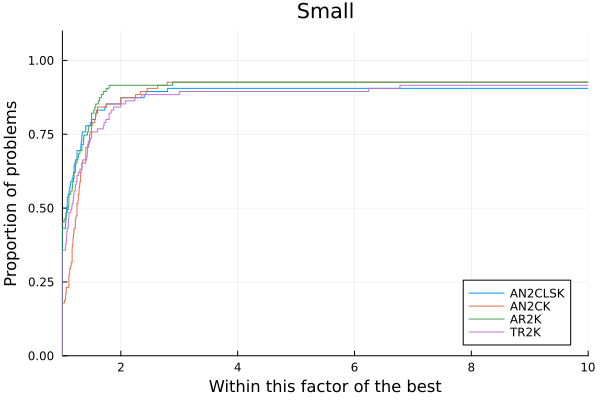}
		\hfil
		\includegraphics[width= 0.5 \linewidth]{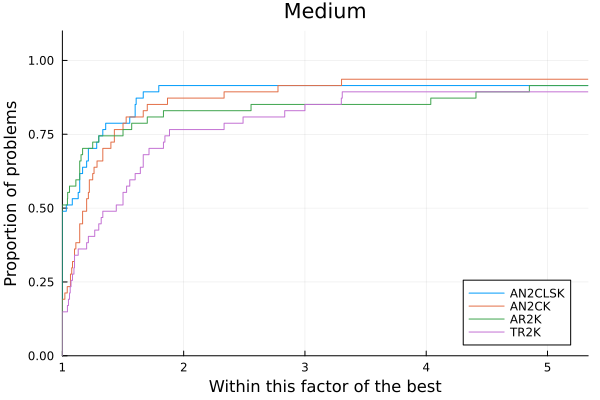}
	}
	\centerline{
		\includegraphics[width= 0.5 \linewidth]{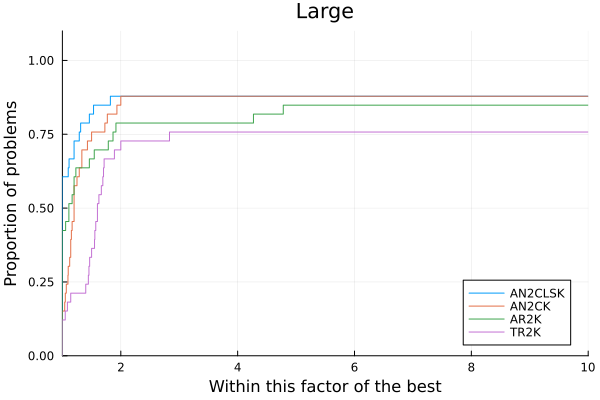}
	} 
	\caption{\label{PlotsK} Performance profile of both \al{AN2CLSK} and \al{AN2CK} on three different sets of problems (small, medium, large)}
\end{figure}

\begin{table}[htb]\footnotesize
	\centering
	\begin{tabular}{|l|r|r|r|r|r|r|}
		\hline
		& \multicolumn{2}{|c|}{small pbs.}& \multicolumn{2}{|c|}{medium pbs.}& \multicolumn{2}{|c|}{largish pbs.}\\
		\hline
		{\tt algo} & $\pi_{\tt algo}$  & $\rho_{\tt algo}$ &
		$\pi_{\tt algo}$ & $\rho_{\tt algo}$ &$\pi_{\tt algo}$ & $\rho_{\tt algo}$ \\
		\hline
		\al{AN2CLSK}  &  0.89 &  91.57  & 0.90 &  91.48 & 0.86 & 87.87 \\
		\al{AN2CK}   &  0.89 &  92.63  & 0.90  &  92.61 & 0.85 & 87.87 \\
		\al{AR2K}   &  0.90 &  92.63  & 0.89  &  91.48 & 0.82 & 84.84 \\
		\al{TR2K}   &  0.88 &  90.57  & 0.84  &  88.23 & 0.76 & 81.81 \\
		\hline
	\end{tabular}
	\caption{\label{tab:statsK} Efficiency and reliability statistics for
		the {\sf{CUTEST}} Julia problems for iterative Krylov subspace solver.}
\end{table}

For the \al{AN2CLSK}  variant, the average ratio of the number of matrix-vector products divided
by the product of the number of iterations and the problem's size (a ratio which is one if every
Lanczos process takes $n$ iterations) is below $0.69$ for small problems, below $0.05$ for medium ones
and below $0.01$ for large ones. Negative curvature directions \eqref{negcurvsteplanczos} are also used, for this variant,
by 0.13\% of the iterations for small problems, 0.02\% of iterations for medium ones and 0\% for
large ones. These statistics do not differ much from those obtained with the \al{AN2CK} variant. 

Again, we see in Table~\ref{tab:statsK} and Figure~\ref{PlotsK} that \al{AN2CLSK} performs on par with the best performing among the other optimization methods 
for the three different sets of problems. 
While these early results are promising and seem to indicate potential for fast second-order methods using local smoothness, the authors are aware that only further experiments
will allow a proper assessment of the method’s true value, both from the number of function/derivatives evaluations and CPU-usage points of view. 

\section{Conclusions and Perspectives}\label{sec:conclusions}

We have proposed a second-order method that can be proved to handle functions whose Hessians only satisfy a local Lipschitz condition. Our algorithm is inspired {by}  \cite{GratJerToin23bIMA} and alternatively uses a regularized Newton step or a negative curvature step when appropriate. 
The proposed algorithm automatically adapts to the problem's geometry and its use does not require prior knowledge of the Lipschitz constant. 
We show that at most $\mathcal{O}(|\log (\epsilon)| \epsilon^{-3/2})$ iterations are required to find an $\epsilon$ first-order approximate point. We have also proposed an algorithmic variant that requires at most $\mathcal{O}\left(\epsilon^{-3}\right)$ iterations to find an approximate second-order critical point.

Two implementations of the framework have been discussed, the first using exact solutions of the involved linear systems and eigenvalue subproblems, the second employing approximations based on nested Krylov subspaces. Numerical experiments suggest that these methods improve slightly on \al{AN2CK} and \al{AN2CE} \cite{GratJerToin23bIMA}, which were already show to be competitive with more standard algorithms using trust-regions or cubic regularization. 

Promising lines of work include inexact or stochastic variants \cite{YaoXuRoosMahoWri22,YaoXuRoosMaho21}, variants where the function value is never evaluated to improve reliability on noisy problems \cite{OFFO-ARp,gratton2024stochastic}, subspace variants  and application to further classes of problems, such as minmax optimization \cite{WangXu24}.

\appendix

\section{Proof of Lemma~\ref{Lipbounds}}
\proof{
	From \cite[Appendix~A.1]{CartGoulToin20b}, we obtain that
	\[
	f(x+s) - f(x) - s^\intercal\nabla_x^1f(x) - \frac{1}{2} s^\intercal\nabla_x^2f(x)s =  \int_{0}^{1} (1-\xi) (s^\intercal\nabla_x^2f(x+\xi s)s-s^\intercal\nabla_x^2f(x)s) \, d\xi
	\]
	Now using basic matrix inequalities, $\|s\|\leq \delta$, $\xi \in
	[0,1]$ and \eqref{LipHessian}, we derive that
	\begin{align*}
		f(x+s) - f(x) - s^\intercal\nabla_x^1f(x) - \frac{1}{2} s^\intercal\nabla_x^2f(x)s &\leq  \int_{0}^{1} (1-\xi) \|\nabla_x^2f(x+\xi s)-\nabla_x^2f(x)\|  \|s\|^2 \, d\xi \\
		&\leq\int_{0}^{1} (1-\xi) (L_0 + L_1 \|\nabla_x^1 f(x)\|) \xi \|s\|^3 \, d\xi \\ 
		&\leq \frac{L_0 + L_1 \|\nabla_x^1 f(x)\|}{6} \|s\|^3
	\end{align*} 
	thus proving \eqref{Lipf}.
	Similarly, the same arguments can be used to prove \eqref{Lipg} since
	\begin{align*}
		\|\nabla_x^1f(x+s) - \nabla_x^1 f(x) - \nabla_x^2 f(x)s\| &= \|\int_{0}^{1} (\nabla_x^2 f(x+\xi s) - \nabla_x^2 f(x))s \, d\xi \| \\
		&\leq \int_{0}^{1} \|\nabla_x^2 f(x+\xi s) - \nabla_x^2 f(x)\|\|s\| \, d\xi \\
		&\leq \int_{0}^1 (L_0 + L_1 \|\nabla_x^1 f(x)\|) \xi \|s\|^2 \, d\xi \\
		&= \frac{(L_0 + L_1 \|\nabla_x^1 f(x)\|) \|s\|^2}{2}.
	\end{align*}
}

\section{Proof of Theorem~\ref{ComplexityboundSecond}}\label{app:so-theory}

As we noted in Section~\ref{so-teaser}, the trial step may be computed in the \al{SOAN2CLS} algorithm using 
Step~2--4 of Algorithm~\ref{AN2Cgen} or \eqref{negcurvstephess}. The 
notations defining the partition of $\ii{k}$ and $\calS_k$ remain relevant,
but we complete them by introducing
\[
\calI^{so} \eqdef \{i\geq 0 \,|\, s_i = s_i^{so} \},
\ms
\calS^{so} \eqdef  \calS \cap \calI^{so},
\ms
\calS_k^{so} \eqdef  \calS_k \cap \calI^{so},
\ms
\calS^{fo} \eqdef  \calS\setminus \calS^{so}
\tim{and}
\calS_k^{fo} \eqdef  \calS_k\setminus \calS_k^{so}.
\]
In addition, for $m\geq \ell \geq 0$, we define
\[
\calS_{\ell,m} \eqdef \calS \cap \iibe{\ell}{m}
\]
and we naturally extend this notation using superscripts identifying
the subsets of $\calS_{\ell,m}$ corresponding to the different  iteration
types identified above.
We  also introduce two index sequences whose purpose is to keep
track of when $s_k= s_k^{fo}$  or $s_k=s_k^{so}$ are used, in the sense that
\[
s_k= s_k^{fo} \tim{for} k \in \bigcup_{i\geq0, p_i\geq0} \iibe{p_i}{q_i-1}
\tim{ and }
s_k= s_k^{so} \tim{for} k \in \bigcup_{i\geq0} \iibe{q_i}{p_{i+1}-1}.
\]
Formally,
\beqn{q0p0def}
p_0 = \left\{\begin{array}{rl}
	0 &\tim{if } \|g_0 \| > \epsilon_1\\
	-1 &\tim{if } \|g_0 \| \le \epsilon_1,
\end{array}\right.
\tim{ and }
q_0 = \left\{\begin{array}{ll}
	\inf \{  k > 0  \mid \|g_{k} \| \leq \epsilon_1 \} &\tim{if } \|g_0 \| > \epsilon_1\\
	0 &\tim{if } \|g_0 \| \le \epsilon_1.
\end{array}\right.
\eeqn
Then
\beqn{qipidef}
p_i \eqdef \inf \{  k > q_{i-1} \mid \|g_k \| > \epsilon_1\}
\tim{ and }
q_i \eqdef \inf \{  k > p_i \mid \|g_k \| \leq \epsilon_1 \}  \tim{ for } i \geq 1.
\eeqn

\noindent
The following lemma states an important decrease property holding when
\req{negcurvstephess} is used. 

\llem{sksoproperties}{
	Suppose that AS.1 and AS.3 hold. Let $k \in \calI^{so}$.
	Then 
	\beqn{hessdecr}
	-g_k^\intercal s_k - \frac {1}{2} s_k^\intercal H_k s_k \geq \frac{1}{2} \sqrt{\sigma_k} |\lambda_{\min}(H_k)| {\|s_k \|^3} = \frac{|\lambda_{\min}(H_k)|}{2 \sigma_{k}}.
	\eeqn
	We also have that if $\|s_k^{so}\| \leq \delta$ 
	\beqn{gkplusonebound}
	\|\nabla_x^1 f(x_k+s_k^{so})\| \leq \frac{L_0+L_1}{2 \sqrt{\sigma_{k} \sigma_{\min}} } +1 + \frac{|\lambda_{\min}(H_k)|}{\sqrt{\sigma_k}}. 
	\eeqn
	
}

\proof{
	First, observe that since $k \in \calI^{so}$, then $\|g_k\| \leq \epsilon_1$. Since termination has not occurred \eqref{SO-termination}, then $\lambda_{\min}(H_k) < - \epsilon_2  < 0$.
	Therefore, we obtain from \eqref{negcurvsecondorder} and the first part of \eqref{negcurvstephess} that
	\[
	g_k^\intercal s_k^{so}+\frac{1}{2}(s_k^{so})^\intercal H_k s_k^{so}
	\leq \frac{1}{2}  \|s_k^{so}  \|^2 u_k^\intercal H_k u_k 
	= \frac{1}{2} \|s_k^{so}  \|^2 \lambda_{\min}(H_k)   
	\leq \frac{\lambda_{\min}(H_k)  }{2} \sqrt{\sigma_{k}} \|s_k^{so}  \|^3,
	\]
	which gives the first inequality in \eqref{hessdecr}. Its second part also follows from the first part of \eqref{negcurvsecondorder}.
	
	Now we turn to the proof of \eqref{gkplusonebound}. Using that $\|s_k^{so}\| \leq \delta$ so that \eqref{Lipg} applies, that $\|g_k\| \leq \epsilon_1 \leq 1$ when $s_k^{so}$ is computed, that $\sigma_{k}\geq \sigma_{\min}$, \eqref{negcurvsecondorder} and \eqref{negcurvstephess}, we derive that
	\begin{align*}
		\|\nabla_x^1 f(x_k+s_k^{so})\| 
		&\leq \frac{L_0+L_1 \|g_k\|}{2} \|s_k^{so}\|^2 + \|g_k + H_ks_k^{so} \| \\
		&\leq \frac{L_0 + L_1}{2 \sigma_{k}} + \|g_k\| + \frac{1}{\sqrt{\sigma_{k}}} \|H_ku_k\| \\
		&\leq \frac{L_0 + L_1}{2 \sigma_{k}} + 1 + \frac{|\lambda_{\min}(H_k)|}{\sqrt{\sigma_{k}}}.
	\end{align*}
	Now using that $\sigma_k \geq \sqrt{\sigma_{\min} \sigma_k}$ gives the second part of the Lemma.	
}
We now prove an analogue of Lemma~\ref{boundsigmak}, now using the negative-curvature 
step as described in \eqref{negcurvsecondorder}-\eqref{negcurvstephess}. We also bound
the sequence of $\|g_{p_i} \|$. 

\llem{hesscurv}{
	Suppose that AS.1, AS.3 and AS.4 hold. Then, for  $k \geq 0$,
	\beqn{sigmamaxhess}
	\sigma_{k} \leq  \max\left(\frac{\kap{max1}}{\epsilon_1}, \frac{\kap{max2}}{\epsilon_2^2} \right),
	\eeqn
	with $\kap{max1} \eqdef \kap{max}$ as defined in \eqref{kappamax} and 
	\beqn{kappamax2}
	\kap{max2} \eqdef \gamma_3 \max\left( \frac{1}{\delta^2}, \frac{(L_0 + L_1)^2}{9(1-\eta_2)^2}\right).
	\eeqn 
}

\proof{
	Note that the results developed when $k \in \calI^{fo}$ in Lemma~\ref{boundsigmak} 
	are still valid with and we now focus on $k \in \calI^{so}$. Supposing that $\sigma_k \geq \frac{1}{\delta^2}$ so that \eqref{Lipf} holds and combining the latter with  \eqref{hessdecr} gives that
	\begin{align*}
		1-\rho_k
		= \frac{f(x_k+s_k)-f(x_k)-g_k^\intercal s_k-\frac{1}{2}s_k^\intercal H_ks_k}
		{-g_k^\intercal s_k - \frac{1}{2} s_k^\intercal H_k s_k} 
		&\leq \frac{( L_0 + L_1)\|s_k^{so}\|^3}{6 (\frac{1}{2}\sqrt{\sigma_k}|\lambda_{\min}(H_k)|\|s_k^{so}\|^3)} \\
		&\leq \frac{L_0+L_1}{3 \sqrt{\sigma_k}|\lambda_{\min}(H_k)|}. 
	\end{align*}
	Thus $\rho_k \geq \eta_2$ provided 
	$\sigma_{k} \geq \frac{(L_0+L_1)^2}{9 \lambda_{\min}(H_k)^2 (1-\eta_2)^2}$. Now injecting that $\sqrt{\sigma_{k}} \geq \frac{(L_0+L_1)}{3 |\lambda_{\min}(H_k)|(1-\eta_2)}$ in \eqref{gkplusonebound} yields
	\[
	\|\nabla_x^1 f(x_k+s_k^{so})\| \leq \frac{3(1-\eta_2) |\lambda_{\min}(H_k)|}{2 \sqrt{\sigma_{\min}}} + 1 + \frac{|\lambda_{\min}(H_k)|}{\sqrt{\sigma_{k}}}.
	\]
	As a consequence, \eqref{nextgradcondhess} must fail and the iteration is successful. The update \eqref{sigmakupdate} therefore ensures that $\sigma_{k+1}\leq \sigma_{k}$ if
	\[
	\sigma_{k} \geq \max\left(\frac{(L_0+L_1)^2}{9 \lambda_{\min}(H_k)^2 (1-\eta_2)} , \frac{1}{\delta^2}, \right).
	\]
	Combining the previous inequality with Lemma~\ref{boundsigmak} and using that $\|g_k\| \geq \epsilon_1$ if $k \in \calI^{fo}$ and $\|\lambda_{\min}(H_k)\| \geq \epsilon_2$ otherwise gives \eqref{sigmamaxhess}. 
}

\noindent
In addition to this lemma, all properties of the different steps
imposed in Section~3 by the mechanism of Algorithm \al{AN2CLS} remain valid. However, \eqref{upperSkdg} in
Lemma~\ref{gdecr} may no longer hold because its proof
relies on the fact that $\|g_k \| \geq \epsilon_1$, which is no longer true.
The purpose of the next lemma is to provide an analogue of
\eqref{upperSkdg} valid for \al{SOAN2CLS}.

\llem{halfnextgradstepsecorder}{
	Suppose that AS.1, AS.3 and AS.4 hold and the \al{SOAN2CLS} algorithm is used. Then
	\begin{align}\label{Skdivfinal}
		|\calS_k^{g\searrow} |
		&\leq \log\left( \frac{\kap{upnewt}}{{\epsilon_1}}\right) \frac{|\calS_k^{decr}|}{\log(2)}  
		+ \log\left(\frac{\kap{upncurv}}{\epsilon_1}\right)  \frac{|\calS_k^{ncurv}|}{\log(2)} \nonumber \\ 
		&\hspace*{5mm}+\left(\frac{|\log(\epsilon_1)| + \log(\kap{gpi})}{\log(2)} +1 \right) 
		(|\calS_k^{so}| + 1)
	\end{align}
	where $\kap{upnewt}$ is defined in Step~0 of the \al{AN2CLS} algorithm,
	$\kap{upncurv}$ defined in \eqref{kapupncurvdef} and 
	$\kappa_{gpi}$ is given by 
	\beqn{kapgpi-def}
	\kap{gpi} 
	\eqdef \max \left[\|g_0\|, \, \frac{3(1-\eta_2)\kappa_B}{2 \sqrt{\sigma_{\min}}} + 1 
	+ \frac{\kappa_B}{\sqrt{\sigma_{\min}}} \right].
	\eeqn
}

\proof{
	We first provide a bound on $\|g_{p_i}\|$ for $p_i \geq 1$. From \eqref{q0p0def} and \eqref{qipidef}, we now know that $p_i - 1 \in \calS^{so}$. Thus \eqref{nextgradcondhess}, the definition of $\kap{k,hess}$ in \eqref{negcurvstephess}, and the facts that $\sigma_{k}\geq \sigma_{\min}$ and that $|\lambda_{\min}(H_k)| \leq \kappa_B$ together ensure that
	\[
	\|g_{p_i}\| \leq \frac{3(1-\eta_2)\kappa_B}{2 \sqrt{\sigma_{\min}}} + 1 + \frac{\kappa_B}{\sqrt{\sigma_{\min}}}.
	\]
	We therefore have that, for all the values $p_i$,
	\beqn{gpibound}
	\|g_{p_i}\| \leq \kap{gpi},
	\eeqn 	
	where $\kap{gpi}$ is defined in \req{kapgpi-def}. We now prove \req{Skdivfinal}.
	If  $\calS_k^{fo}$ is empty, its subset $\calS_k^{g\searrow}$ is also empty and
	\req{Skdivfinal} trivially holds.
	If $\calS_k^{fo}$ is not empty, we see from the definitions
	\eqref{q0p0def}-\eqref{qipidef} that, for some $m\geq0$ depending on $k$,
	\beqn{0kdivbiggerepsilon}
	\calS_k^{fo}
	= \calS_{k} \cap \{ \|g_k \| > \epsilon_1 \}
	= \left(\bigcup_{i=0, p_i\geq0}^{m-1} \iibe{p_i}{q_i -1}\right) \cup \iibe{p_{m}}{k}.
	\eeqn
	Note that the last set in this union is empty unless $k\in\calS^{fo}$, in which case $p_m \geq 0$. 
	Suppose first that the set of indices corresponding to the union in
	brackets is non-empty and let $i$ be an index in this set.  Moreover,
	suppose also that $p_i < q_i-1$. Using \req{gpibound} and the facts that
	$\| g_{q_{i} -1} \| > \epsilon_1$, that the gradient only changes at successful iterations and that 
	$\calS_{p_i,q_i-2}=\calS_{p_i,q_i-2}^{ncurv}\cup\calS_{p_i,q_i-2}^{g\searrow}\cup\calS_{p_i,q_i-2}^{decr}$ the fact that  \eqref{gkplusoncurvbound} holds for $i \in \calS_{p_i,q_i-2}^{ncurv}$, and that \eqref{nextgradcond} fails for $i \in \calS_{p_i,q_i-2}^{decr} \cup \calS_{p_i,q_i-2}^{ncurv}$ with $\epsilon = \epsilon_1$ in both cases, we  derive that 
	\begin{align*}
		\frac{\epsilon_1}{\kappa_{gpi}} 
		&\leq \frac{\|g_{q_i -1 } \|}{\| g_{p_i}\|} 
		= \prod_{j = p_i}^{q_i - 2} \frac{\|g_{j+1} \|}{\|g_j \|} 
		= \prod_{j \in \calS_{{p_i, q_i - 2}} } \frac{\|g_{j+1} \|}{\|g_j \|}\\  
		&= \prod_{j \in \calS_{p_i, q_i - 2}^{decr}}  \frac{\|g_{j+1}\|}{\|g_j \|}
		\prod_{j \in \calS_{p_i, q_i - 2}^{ncurv}}  \frac{\|g_{j+1} \|}{\|g_j \|} 
		\prod_{j \in \calS_{p_i, q_i - 2}^{g\searrow}} \frac{\|g_{j+1} \|}{\|g_j \|}\\
		&\leq \left(\frac{\kap{upnewt}}{\epsilon_1}\right)^{|\calS_{p_i, q_i-2}^{decr} |} 
		\times  \frac{1}{2^{|\calS_{p_i, q_i-2}^{g\searrow} |}} 
		\times  \left(\frac{\kap{upncurv}} {\epsilon_1}\right)^{|\calS_{p_i,q_i-2}^{ncurv} |}.
	\end{align*}
	Rearranging terms, taking the log, using the inequality 
	$|\calS_{p_i,q_i-2}^{g\searrow}| \geq |\calS_{p_i,q_i-1}^{g\searrow}|-1$ and
	dividing by $\log(2)$ then gives that
	\[
	(|\calS_{p_i,q_i-1}^{g\searrow}|-1) + \frac{\log (\epsilon_1)-\log(\kappa_{gpi})}{\log(2)}
	\leq \frac{\log\left(\frac{\kap{upnewt}}{\epsilon_1}\right)}{\log(2)} |\calS_{p_i, q_i - 2}^{decr}| 
	+ \frac{\log\left(\frac{\kap{upncurv}} {\epsilon_1}\right)}{\log(2)} |\calS_{p_i, q_i - 2}^{ncurv}|.
	\]
	Further rearranging and using the fact that
	$|\calS_{p_i, q_i - 2}| \leq |\calS_{p_i, q_i - 1}|$ for the different
	types of step, we obtain that
	\beqn{Skdivpiqi}
	|\calS_{p_i, q_i - 1}^{g\searrow}|
	\leq \frac{\log\left(\frac{\kap{upnewt}}{\epsilon_1}\right)}{\log(2)} |\calS_{p_i,q_i-1}^{decr}|
	+ \frac{\log\left(\frac{\kap{upncurv}} {\epsilon_1}\right)}{\log(2)}|\calS_{p_i,q_i-1}^{ncurv}|
	+ \frac{|\log (\epsilon_1)| + \log(\kappa_{gpi}) }{\log(2)}  + 1.
	\eeqn
	If now $p_i = q_i-1$, then clearly $|\calS_{p_i,q_i-1}^{g\searrow}|\leq 1$ and \req{Skdivpiqi} also holds.
	Using the same reasoning  when $\iibe{p_m}{k}$ is non-empty, we  derive that,
	\beqn{Skdivpmk}
	|\calS_{p_m, k}^{g\searrow}|
	\leq \frac{\log\left(\frac{\kap{upnewt}}{\epsilon_1}\right)}{\log(2)} |\calS_{p_m,k}^{decr}|
	+ \frac{\log\left(\frac{\kap{upncurv}} {\epsilon_1}\right)}{\log(2)}|\calS_{p_m,k}^{ncurv}|
	+ \frac{|\log (\epsilon_1)| + \log(\kappa_{gpi}) }{\log(2)}  + 1,
	\eeqn
	and this inequality also holds if  $\iibe{p_m}{k}=\emptyset$ since 
	$\calS_{p_m, k}^{g\searrow}\subseteq\iibe{p_m}{k}$.
	Adding now \eqref{Skdivpiqi} for $i \in \iiz{m}$ and \req{Skdivpmk} to take 
	\req{0kdivbiggerepsilon} into account gives that
	\begin{align*}
		|\calS_{ k}^{g\searrow}|
		\leq \frac{\log\left(\frac{\kap{upnewt}}{\epsilon_1}\right)}{\log(2)} |\calS_{ k}^{decr}|
		+ \frac{\log\left(\frac{\kap{upncurv}} {\epsilon_1}\right)}{\log(2)} |\calS_{k}^{ncurv}| 
		+ \left(\frac{|\log (\epsilon_1)| + \log(\kappa_{gpi}) }{\log(2)}+1\right) (m+1).
	\end{align*}
	Because \eqref{0kdivbiggerepsilon} divides $\calS_{k}^{fo}$
	into $m+1$ consecutive sequences, these sequences are then separated
	by at least one second-order step, so that $m \leq |\calS_{ k}^{so}|$
	and \req{Skdivfinal} follows. 
}

\noindent
Equipped with this last lemma and the results of
Sections~\ref{thealgo-s} and \ref{complexity-s}, we
may finally establish the worst-case complexity
of the \al{SOAN2CLS} algorithm and prove Theorem~\ref{ComplexityboundSecond} itself.
\vspace*{3mm}
\proof{
	Since the  bounds \eqref{Skcurvbound} and \eqref{Skdecrbound} depended on the value of 
	$\sigma_{\max}$ obtained in \eqref{sigmamax}, they must be rederived since a new bound 
	\eqref{sigmamaxhess} holds. 
	
	We start by providing a bound on $|\calS_{k}^{ncurv}|$. 
	Using AS.2, \eqref{quaddecrcurv}, the fact that $\|g_i\| \geq \epsilon_1$ before termination and \eqref{sigmamaxhess}, we deduce that,
	\begin{align*}
		f(x_0) - f_{\rm low}
		\geq \sum_{i \in \calS_k} f(x_i) - f(x_{i+1}) &\geq  \sum_{i \in \calS_k^{ncurv}} f(x_i) - f(x_{i+1}) \\
		&\geq \sum_{i \in \calS_k^{ncurv}}  \, \frac{\eta_1 \theta^3 \kappa_C^3 }{2 \sqrt{\sigma_{i}}}\,\|g_i\| 
		\geq \frac{\eta_1 \theta^3 \kappa_C^3 \epsilon_1 |\calS_k^{ncurv}| }{2 \max\left(\sqrt{\frac{\kap{max1}}{\epsilon_1}}, \, \frac{\sqrt{\kap{max2}}}{\epsilon_2}\right)} ,
	\end{align*}
	and hence after rearranging and using Young's inequality with $p = \frac{3}{2}$ and $q=3$, we derive that
	\beqn{Skcurvboundhess}
	|\calS_k^{ncurv}| 
	\leq \frac{2 (f(x_0) - f_{\rm low})  \max\left[\sqrt{\kap{max1}} \epsilon_1^{-3/2}, \sqrt{\kap{max2}} (\frac{2}{3}\epsilon_1^{-3/2} + \frac{1}{3} \epsilon_2^{-3}) \right]}{\eta_1 \theta^3 \kappa_C^3} 
	\leq \kap{ncurvhess} (\epsilon_1^{-3/2} + \epsilon_2^{-3})
	\eeqn
	with $\kap{ncurvhess}$ defined in \eqref{kaphessdef}.
	
	We now provide a bound on $|\calS_{ k}^{decr}|$.
	Using AS.2, \eqref{quaddecrneig}, that $\|s_i^{decr}\| \geq \frac{1}{\sqrt{\sigma_{k}} \kap{slow}}$, the fact that $\|g_i\|\geq \epsilon_1$ before termination, and  \eqref{sigmamaxhess} yields that
	\begin{align*}
		f(x_0) - f_{\rm low}	&\geq \sum_{i \in  \calS_k^{decr}} f(x_i) - f(x_{i+1})
		\geq \sum_{i \in  \calS_k^{decr}}\eta_1  \, \sqrt{\sigma_i} \|g_i\| \|s_i\|^2 \\ &\geq \sum_{i \in \calS_k^{decr}}\eta_1     \frac{\|g_i\| }{\sqrt{\sigma_{i}} \kap{slow}^2} 
		\geq  \frac{\epsilon_1 |\calS_k^{decr}|\eta_1 }{\max\left(\sqrt{\frac{\kap{max1}}{\epsilon_1}}, \, \frac{\sqrt{\kap{max2}}}{\epsilon_2}\right) \kap{slow}^2}. 
	\end{align*}
	Again, by the same arguments used to prove \eqref{Skcurvboundhess}, we derive that
	\beqn{Skdecrboundhess}
	|\calS_k^{decr}| \leq \frac{ (f(x_0) - f_{\rm low}) \kap{slow}^2 \max\left(\sqrt{\kap{max1}} \epsilon_1^{-3/2}, \sqrt{\kap{max2}} (\frac{2}{3}\epsilon_1^{-3/2} + \frac{1}{3} \epsilon_2^{-3}) \right)}{\eta_1 } \leq \kap{decrhess} (\epsilon_1^{-3/2} + \epsilon_2^{-3})
	\eeqn
	with $\kap{decrhess}$ defined in \eqref{kaphessdef}.
	
	We finally prove a bound on $|\calS_k^{so}|$.
	Using AS.2 and the lower bound on the decrease of the function values
	\eqref{hessdecr}, \eqref{sigmamaxhess} and the fact that $|\lambda_{\min}(H_i)| \geq \epsilon_2$, we conclude that
	\[
	f(x_0) - f_{\rm low}
	\geq  \sum_{i \in \calS_k^{so}} f(x_i) - f(x_{i+1})
	\geq \sum_{i \in \calS_k^{so}} \, \frac{\eta_1 |\lambda_{\min}(H_i)|}{2 \sigma_{i}}\, \geq |\calS_k^{so}| \, \frac{\eta_1  \epsilon_2 }{2\max \left(\frac{\kap{max1}}{\epsilon_1}, \frac{\kap{max2}}{\epsilon_2^{2}}\right)}.
	\]
	Rearranging the last inequality and using Young's inequality, we obtain that
	\beqn{Skhessbound}
	|\calS_k^{so}| \leq  \frac{ 2 (f(x_0) - f_{\rm low})  \max\left( {\kap{max2}} \epsilon_2^{-3}, \kap{max1} (\frac{2}{3}\epsilon_1^{-3/2} + \frac{1}{3} \epsilon_2^{-3}) \right)}{\eta_1 } \leq \kap{so} (\epsilon_1^{-3/2} + \epsilon_2^{-3})
	\eeqn
	where $\kap{so}$ is defined at \eqref{kappastarlogkapso}.
	
	Substituting now \eqref{Skhessbound}, \eqref{Skdecrboundhess} and
	\eqref{Skcurvboundhess} in the bound \eqref{Skdivfinal} on $\calS_k^{g\searrow}$ yields
	\begin{align*}
		|\calS_{k}^{g\searrow}| &\leq
		\log\left(\frac{\kap{upnewt}}{{\epsilon_1}}\right) \frac{\kap{decrhess} (\epsilon_1^{-3/2} + \epsilon_2^{-3})}{\log(2)}  + 
		\log\left(\frac{\kap{upncurv}}{{\epsilon_1}}\right) \frac{\kap{curvhess} (\epsilon_1^{-3/2} + \epsilon_2^{-3})}{\log(2)} \nonumber \\ &+ \left(\frac{|\log(\epsilon_1)| + \log(\kap{gpi})}{\log(2)} + 1 \right) \Big(\kap{so}(\epsilon_1^{-3/2} + \epsilon_2^{-3}) + 1\Big).
	\end{align*}
	Combining this last inequality with  \eqref{Skhessbound}, \eqref{Skdecrboundhess} 
	and \eqref{Skcurvboundhess} in
	$|\calS_k| = |\calS_k^{g\searrow}| + |\calS_k^{ncurv}| + |\calS_k^{so}| + |\calS_k^{decr}|$
	and using the definitions of $\kappa_{\star,hess}$  \eqref{kappastarhess} and $\kappa_{\star,log}$ \eqref{kappastarlogkapso}  gives that
	\[
	|\calS_{ k} | \leq \kappa_{\star,hess} \left(\epsilon_1^{-3/2} + \epsilon_2^{-3}\right)  + \kappa_{\star,log} |\log(\epsilon_1)| (\epsilon_1^{-3/2} + \epsilon_2^{-3}) + \frac{|\log(\epsilon_1)| + \log(\kap{gpi})}{\log(2)} + 1,
	\]
	which proves the first part of the theorem. 
	The second part follows from the last inequality, Lemma~\ref{SvsU} and the bound \eqref{sigmamaxhess}.
}





\end{document}